\documentclass[UTF8,12pt]{article}
\usepackage{upgreek}
\usepackage{listings}
\usepackage{amsmath,amssymb}
\usepackage{longtable}
\usepackage{amsthm}
\usepackage{multirow}
\usepackage{mathrsfs}
\usepackage{pgf,tikz}
\usepackage{pgfplots}
\usepackage{makecell}
\usepackage{graphicx}
\usepackage{subfigure}
\usepackage{caption}
\usepackage{algorithm}
\usepackage[noend]{algorithmic}
\usepackage{bm}
\usepackage{bbm}
\usepackage{lscape}
\usetikzlibrary{arrows}
\usepackage{geometry}
\usepackage{booktabs}
\usepackage{url} 
\usepackage{ragged2e}
\usepackage{hyperref}
\usepackage{float}
\renewcommand\leq{\leqslant}
\renewcommand\geq{\geqslant}
\renewcommand\pi{\uppi}

\newtheorem{theorem}{Theorem}[section]

\newtheorem{lemma}[theorem]{Lemma}
\newtheorem{corollary}[theorem]{Corollary}
\newtheorem{remark}[theorem]{Remark}
\newtheorem{definition}[theorem]{Definition}

\usepackage{array}
\newcolumntype{L}[1]{>{\RaggedRight\arraybackslash}p{#1}} 
\newcolumntype{C}[1]{>{\Centering\arraybackslash}p{#1}} 
\newcolumntype{R}[1]{>{\RaggedLeft\arraybackslash}p{#1}} 

\newcolumntype{N}[1]{>{\RaggedLeft\arraybackslash $ }p{#1}<{$}} 
\newcolumntype{M}{>{$}r<{$}} 
\newcolumntype{P}{>{$}l<{$}} 

\makeatletter
\@addtoreset{equation}{section}
\makeatother

\usepackage{xcolor}

\usepackage{ulem}

\begin{document}
\title{Application of Causal Inference Techniques to the Maximum Weight Independent Set Problem}
\author{Jianfeng Liu \footnotemark[2] \footnotemark[4] \and Sihong Shao\footnotemark[3] $\thanks{Author to whom correspondence should be addressed: sihong@math.pku.edu.cn}$ \and Chaorui Zhang\footnotemark[4]}

\renewcommand{\thefootnote}{\fnsymbol{footnote}}
\footnotetext[2]{Department of Mathematical Sciences, Tsinghua University, Beijing 100084, P.R. China}
\footnotetext[3]{CAPT, LMAM and School of Mathematical Sciences, Peking University, Beijing 100871, P.R. China}
\footnotetext[4]{Theory Lab,~Central Research Institute,~2012 Labs, Huawei Technologies Co., Ltd.}
\date{\today}
\maketitle

\begin{abstract}
A powerful technique for solving combinatorial optimization problems is to reduce the search space without compromising the solution quality by exploring intrinsic mathematical properties of the problems. For the maximum weight independent set~(MWIS) problem, using an upper bound lemma which says the weight of any independent set not contained in the MWIS is bounded from above by the weight of the intersection of its closed neighbor set and the MWIS,  we give two extension theorems --- independent set extension theorem and vertex cover extension theorem. With them at our disposal, two types of causal inference techniques~(CITs) are proposed on the assumption that a vertex is strongly reducible~(included or not included in all MWISs) or reducible~(contained or not contained in a MWIS). One is a strongly reducible state-preserving technique, which extends a strongly reducible vertex into a vertex set where all vertices have the same strong reducibility. The other, as a reducible state-preserving technique, extends a reducible vertex into a vertex set with the same reducibility as that vertex and creates some weighted packing constraints to narrow the search space. Numerical experiments show that our CITs can help reduction algorithms find much smaller remaining graphs, improve the ability of exact algorithms to find the optimal solutions and help heuristic algorithms produce approximate solutions of better quality. In particular, detailed tests on $12$ representative  graphs generated from datasets in Network Data Repository demonstrate that, compared to the state-of-the-art algorithms, the size of remaining graphs is further reduced by more than $32.6\%$, and the number of solvable instances is increased from $1$ to $5$. 

\vspace*{4mm}
\noindent {\bf AMS subject classifications:}
05C69; 
68W40; 
90C06; 
90C27; 
90C57 

\noindent {\bf Keywords:}
maximum weight independent set; independent set extension; vertex cover extension; causal inference techniques; reduction algorithm; exact algorithm; heuristic algorithm; Network Data Repository.
\end{abstract}
\section{Introduction}
\label{sec:Introduction}
Let $G = (V, E, w)$ be an undirected vertex-weighted graph, where each vertex $v \in V$ is associated with a weight $w(v) \in \mathbb{R}^{+}$. A subset $I \subseteq V$ is called an independent set if its vertices are pairwise non-adjacent, and the vertex cover of graph $G$ is a subset of vertices $VC \subseteq V$ such that every edge $e \in E$ is incident to at least one vertex in subset $VC$. Independent set and vertex cover are two complementary concepts in graph and can be transformed into each other on demand~\cite{xu2016new}. The maximum weight independent set~(MWIS) problem is to find the independent set of largest weight among all possible independent sets and the weight of a MWIS of graph $G$ is denoted by $\alpha_w(G)$, while the minimum weight vertex cover~(MWVC) problem asks for the vertex cover with the minimum weight. Furthermore, if subset $I \subseteq V$ is a MWIS, then subset $VC = V \backslash I$ is a MWVC, and vice versa~\cite{cai2018improving,xu2016new}. The MWIS problem is an extension of the maximum independent set~(MIS) problem, which is a classic NP-hard problem~\cite{fomin2009measure,coniglio2022optimizing}. It can be applied to various real-world problems, such as information retrieval~\cite{balas1986finding}, computer vision~\cite{feo1994greedy}, combinatorial auction problem~\cite{xu2016new} and dynamic map labeling problem~\cite{liao2016approximation}. Due to its wide range of practical applications, the research on efficient algorithms for computing the MWIS is of great significance. Most previous work are focused on heuristic algorithms to find near-optimal solutions in reasonable time~\cite{shyu2004ant,pullan2009optimisation,cai2018improving,nogueira2018hybrid}, while exact algorithms, usually referring to Branch-and-Bound~(B\&B) methods~\cite{babel1994fast,warren2006combinatorial, avenali2007resolution, san2019new}, become infeasible when the size of problem increases. 

Recently, it has been well demonstrated that reduction rules~(a.k.a. kernelization) are very effective in practice for solving the MIS problem~\cite{strash2016power}. These rules mine the structural properties of underlying graph and reduce the search space by such as removing vertices, contracting subgraphs, restricting the set of independent sets, etc., to produce a smaller kernel graph such that the MIS of the original graph can be recovered from the MIS of the kernel. After integrating them,  some state-of-the-art exact solvers are able to solve the MIS problem on many large real networks~\cite{dzulfikar2019pace}. These solvers can be usually divided into two types: One performs the kernelization only once and runs the B\&B algorithm~\cite{sewell1998branch,li2018incremental} on the kernelized instance, while the other joins hands with the Branch-and-Reduce~(B\&R) algorithm~\cite{plachetta2021sat} and performs reduction in every branch of the search tree. As for those instances that can’t be solved exactly, high-quality solutions can be found by combining kernelization with local search~\cite{chang2017computing,dahlum2016accelerating}. Moreover, when a vertex is selected for branching in the branching process of the B\&R algorithm, if it is assumed to be in all MISs, then its satellite set will also be in all MISs~\cite{kneis2009fine}, while its mirror set will be removed directly from the graph, if it is assumed not to be in all MISs~\cite{fomin2009measure}. Further, a conflict analysis on the assumption that a vertex is in all MISs can be also plugged in to find some contradictions and the concept of ``unconfined/confined vertices" was introduced~\cite{xiao2013confining}. Later, an auxiliary constraint called packing constraint was proposed  to accelerate the B\&R algorithm by simply exploring branches that satisfy all packing constraints~\cite{akiba2016branch}. The central idea behind all these attempts for the MIS problem involves a state-preserving technique which starts from a vertex, named the starting vertex for convenience, and then finds a vertex set with the same state as the starting vertex to reduce the search space, thereby implying that some subsequent operations can be implemented on the resulting vertex set instead of only on the starting vertex. For the MWIS problem, similar state-preserving techniques are rarely used except for a recent work using unconfined/confined vertices~\cite{xiao2021efficient}, though some simple and fast reduction rules have been used in  B\&R algorithms~\cite{lamm2019exactly,xiao2021efficient}. To this end, we devote ourselves into developing state-preserving techniques for the MWIS problem in this work. The state of the starting vertex we consider can be
\begin{itemize}
\item strongly reducible, meaning that the vertex is included in all MWISs/MWVCs; or 
\item reducible, meaning that the vertex is contained in a MWIS/MWVC.
\end{itemize}
Considering that the assumed state of the starting vertex must be used to analyze its local structure to obtain inference results, these targeted state-preserving techniques are called causal inference techniques~(CITs).  Inspired by their success in solving the MIS problem, we will systematically develop CITs to solve the MWIS problem by analyzing intrinsic mathematical properties of underlying graph. More specifically, our main contributions are in three aspects as follows. 

First, by virtue of the upper bound lemma, i.e., the weight of any independent set not contained in the MWIS is bounded from above by the weight of the intersection of its closed neighbor set with the MWIS, two extension theorems are developed. With them, we propose a series of CITs which have been rarely used previously in the MWIS problem. According to the state of the starting vertex, our CITs can be divided into two categories. The first type is a strongly reducible state-preserving technique. We first assume that the starting vertex is strongly reducible, and then try to extend this vertex to obtain a vertex set with the same strong reducibility. If the upper bound lemma is not satisfied in this process, then this contradicts the assumption, and the starting vertex can be removed from the graph directly. Otherwise, combined with the state-preserving result obtained from the previous process, we continue to search for a set called the simultaneous set, which is either included in a MWIS or contained in a MWVC. The second type is a reducible state-preserving technique. Under the assumption that the starting vertex is reducible, a vertex set with the same reducibility can be obtained by extending from this vertex. Moreover, if this vertex is selected for branching in the B\&R algorithm, with the upper bound lemma, an inequality constraint called weight packing constraint will be created to restrict subsequent searches.
     
Next, according to the characteristics of the proposed CITs, we integrate them into the existing algorithmic framework. The first type of CIT can be used to  design reduction rules to simplify graph. These reduction rules are integrated into the existing reduction algorithm. In the B\&R algorithm, when a vertex is selected to branch, a vertex set and a weight packing constraint depending on the assumed state of the vertex can be obtained from state-preserving results of two types of CITs. The vertex set is used to further simplify the corresponding branch, while we can prune branches that violate constraints and simplify the graph by maintaining all created weight packing constraints. During the local search process of the heuristic algorithm, when the state of a vertex needs to be changed, all vertex states in the vertex set obtained by the second type of CIT will also be modified to be the same as that vertex, which expands the area of local search and improves the ability of local search to find better local optima. 

Numerical experiments on $12$ representative graphs generated from datasets in Network Data Repository show that the performance of various algorithms is greatly improved after integrating our CITs. The size of the kernel obtained by the resulting reduction algorithm is greatly reduced. In addition, compared to the state-of-the-art exact algorithm, the number of solvable instances have been increased from $1$ to $5$. And the ability of the heuristic algorithm to find better local optimal solutions is significantly improved. These experimental results form the third major contribution of this paper.

\begin{table}[ht]
    \centering
    \resizebox{\textwidth}{25mm}{
    \begin{tabular}{|c|c||c|c|} 
    \hline
    $G = (V, E, w)$                                                                  & \multicolumn{3}{c|}{an undirected vertex-weight  graph $G$ with vertex set $V$, edge set $E$ and vertex weight function $w: V \rightarrow \mathbb{R}^{+}$ } \\ \hline
    $N(v) = \{ u \in V \vert \{ u, v \} \in E \}$ & the neighbor set of vertex $v$  &  $N[v] = N(v) \cup \{ v\}$         &  the closed neighbor set of vertex $v$ \\ \hline
    $N(S) = (\bigcup\limits_{v \in S} N(v)) \backslash S$ & the open neighbor set of set $S$   & $N[S] = N(S) \cup \{ S \}$ &  the closed neighbor set of  set $S$  \\ \hline
    $\vert S \vert$ &  the size of set $S$  &  $w(S) = \sum\limits_{v \in S} w(v)$                        &  the weight of all vertices in set $S$   \\ \hline
     $d(v)$ & the degree of a vertex $v$            & dist$(u, v)$  &  \makecell[c]{ the minimum number of edges in the path from vertex $u$ to vertex $v$}  \\ \hline  
     \makecell[c]{ $N^{l}(v) = \{ u \vert$ dist$(u,v) = l \}$}              & \makecell[c]{the set of vertices at distance $l$ from vertex $v$, \\in particular, $N^{1}(v) = N(v)$} &
        \makecell[c]{$G[S] = (S, E_S, w)$, \\$ \forall e =  \{ u, v \} \in E_S$, $u, v  \in S$} &\makecell[c]{the subgraph induced by a non-empty vertex subset $S$ of $V$}                           \\ \hline
        $\alpha(G)$ & the size of a MIS of unweight graph $G$ & $\alpha_w(G)$   &                             the weight of a MWIS of  graph $G$ \\ \hline
        $A_I$ & the set of all MWISs in graph $G$   &$A_C$                                            &   the set of all MWVCs in graph $G$  \\ \hline
    $S \blacktriangleleft A_I$& set $S$ is an independent set and is included in all MWISs &  $C \vartriangleleft  A_C$ &  set $C$ is contained in all MWVCs              \\ \hline
    vertex $v$ is strongly reducible & vertex $v$ is included in all MWISs/MWVCs &  vertex $v$ is reducible & vertex $v$ is contained in a MWIS/MWVC          \\ \hline 
    
    vertex $v$ is strongly inclusive & vertex $v$ is included in all MWISs     & vertex $v$ is strongly sheathed & vertex $v$ is contained in all MWVCs     \\ \hline 
    vertex $v$ is inclusive & vertex $v$ is included in a MWIS     & vertex $v$ is sheathed & vertex $v$ is contained in a MWVC     \\ \hline 
   
    set $S$ is strongly inclusive  & set $S$ is an independent set and is included in all MWISs     & set $S$ is strongly sheathed & set $S$ is contained in all MWVCs     \\ \hline 
    set $S$ is inclusive & set $S$ is an independent set and is included in a MWIS    & set $S$ is sheathed & set $S$ is contained in a MWVC     \\ \hline 
    independent set $S$ is strongly  exclusive  & independent set $S$ is not contained in all MWIS      &   independent set $S$ is exclusive  & independent set $S$ is not contained in a MWIS \\ \hline 
    a set $S$ called a simultaneous set & \multicolumn{3}{c|}{ set $S$ is either included in a MWIS or contained in a MWVC} \\ \hline 
    \end{tabular}}
    \caption{Notations used throughout the paper.}
    \label{tab: Relevant notation used throughout the paper}
\end{table}

Relevant notations used in this work are given in Table~\ref{tab: Relevant notation used throughout the paper}
and the rest of the paper is organized as follows. We present two extension theorems in Section~\ref{sec:Two Extension Theorems}
and detail CITs in Section~\ref{sec:Causal Inference Technique}. 
How the CITs are combined with existing algorithmic frameworks is described in Section~\ref{sec:Integrate CITs into Existing Algorithmic Frameworks}. Extensive numerical tests are carried out in Section~\ref{sec:Experiments} to verify the performance improvement of integrating our CITs into existing algorithmic frameworks in terms of efficiency and accuracy. 
The paper is concluded in Section~\ref{sec:Conclusion} 
with a few remarks. 

\section{Two Extension Theorems}
\label{sec:Two Extension Theorems}
The theoretical cornerstones of CITs in this paper are two extension theorems: independent set extension theorem and vertex cover extension theorem. Before delineating them, we need to have a deep understanding of the local structure of the MWIS and first give the upper bound lemma.

\begin{lemma}[upper bound lemma]
\label{th:Local Structure of Independent Set}
Let set $I_C$ be an independent set in the graph.
\begin{itemize}
\item [$(a)$] Suppose there is an $I_w \in A_I$ such that $I_C \not \subseteq  I_w$, then $w(I_w \cap N[I_C]) \geq w(I_C)$ holds.
\item [$(b)$] Assume that $I_C \not \subseteq I, \forall I \in A_I$ holds, then it satisfies: $w(I_C) < w(I \cap N[I_C]), \forall I \in A_I$. 
\end{itemize}
\end{lemma}
\proof{Proof}
We first prove $(a)$ by contradiction. If not, we can obtain an independent set $ I_w^{\prime} = (I_w  \backslash N[I_C]) \cup (I_C) $ such that $w(I_w^{\prime}) = w(I_w) + w(I_C) - w(I_w \cap N[I_C]) > w(I_w)$, a contradiction. 

Next, we consider $(b)$. If there is an $I_1 \in A_I$ such that $w(I_C) \geq w(I_1 \cap N[I_C]))$, we can construct an independent set $I^{\prime}_1 = (I_1  \backslash N[I_C]) \cup I_C$ satisfying $w(I^{\prime}_1) = w(I_1) + w(I_C) - w(I_1 \cap N[I_C]) \geq w(I_1)$. Then $I^{\prime}_1 \in A_I$ and $I_C \subseteq  I^{\prime}_1$, which leads to a contradiction.
\endproof

The upper bound lemma describes such a property: For any independent set that is~(strongly) exclusive, the weight of the intersection of its closed neighbor set with the MWIS is the upper bound on its weight. With it, the independent set extension theorem can be introduced as follows.
\begin{theorem}[Independent Set Extension Theorem]
\label{th:independent set extension theorem}
Let sets $IS$ and $S$ be two independent sets in the graph.
\begin{itemize}
\item [$(a)$] Assume that there exists an $I_w \in A_I$ such that $IS \subseteq I_w$. If there is an independent set $IS^{\prime} \subseteq N(IS)$ such that $w(IS^{\prime}) > w(IS \cap N(IS^{\prime}))$, then there exists an independent set $IS^{\prime \prime} \subseteq N(IS^{\prime}) \backslash N[IS]$ satisfying the inequality: $w(IS^{\prime}) \leq w(IS \cap N(IS^{\prime})) + w(IS^{\prime \prime})$. In addition, $IS \cup IS^{\prime \prime} \subseteq I_w$ if such $IS^{\prime \prime}$ is unique.
\item [$(b)$] Suppose $S \blacktriangleleft A_I$, then for any independent set $S^{\prime} \subseteq N(S)$, there is an independent set $S^{\prime \prime} \subseteq N(S^{\prime}) \backslash N[S]$ such that $w(S^{\prime}) < w(S \cap N(S^{\prime})) + w(S^{\prime \prime})$. Besides, if such $S^{\prime \prime}$ is unique, then $S \cup S^{\prime \prime} \blacktriangleleft A_I$.
\end{itemize}
\end{theorem}
\proof{Proof}
We first consider the proof of $(a)$, and it is obvious that $IS^{\prime} \not \subseteq I_w$.
In view of the fact that the relationship between $I_w$ and $N[IS^{\prime}]$ satisfies: $I_w \cap N[IS^{\prime}] = I_w \cap N(IS^{\prime})  = (IS \cap N(IS^{\prime})) \cup (I_w \cap (N(IS^{\prime}) \backslash N[IS]))$ and by the upper bound lemma, we can get: 
$ w(IS \cap N(IS^{\prime})) + w(I_w \cap (N(IS^{\prime}) \backslash N[IS])) = w(I_w \cap N(IS^{\prime})) = w(I_w \cap N[IS^{\prime}])\geq w(IS^{\prime}).$
Thus, the existence of such $IS^{\prime \prime}$ is proved. Furthermore, assuming that such $IS^{\prime \prime}$ is unique, then $IS^{\prime \prime} = I_w \cap (N(IS^{\prime}) \backslash N[IS])$ and $IS \cup IS^{\prime \prime} \subseteq I_w$. 

Similar ideas can be used to prove $(b)$. Obviously $S^{\prime} \not \subseteq I, \forall I \in A_I$ holds, so from the upper bound lemma, it can be directly obtained: $\forall I \in A_I,  w(I \cap N[S^{\prime}])  > w(S^{\prime})$. 
Further, by considering that the relationship between $I$ and $N[S^{\prime}]$ satisfies: $I \cap N[S^{\prime}] = I \cap N(S^{\prime}) = (S \cap$ $N(S^{\prime}))  \cup (I \cap (N(S^{\prime}) \backslash N[S]))$,  we prove the existence of such $S^{\prime \prime}$. Also, if such $S^{\prime \prime}$ is unique, the following result holds: $S^{\prime \prime} = I \cap (N(S^{\prime}) \backslash N[S]), \forall I \in A_I$, and then $S \cup S^{\prime \prime} \blacktriangleleft A_I$. 
\endproof

The independent set extension theorem gives a method for extending independent set that is~(strongly) inclusive: Try to find an independent set to add to the extended independent set, and that independent set is the only one that guarantees that the upper bound lemma is satisfied in the local structure of the extended independent set. Next, with the help of the upper bound lemma,  the vertex cover extension theorem is given below.
\begin{theorem}[Vertex Cover Extension Theorem]
\label{th:vertex cover extension theorem}
Let sets $IC$ and $C$ be two vertex subsets in the graph.
\begin{itemize}
    \item[$(a)$] Suppose set $IC \subseteq VC_w$, then the vertices in $IC$ have the property: $\forall p \in IC$, $w(p) \leq \alpha_w(G[N(p)\backslash IC])$. Also, for a vertex $v\in IC$ and a vertex $u \in N^2(v)$, $IC \cup \{ u\}  \subseteq VC_w$ holds if the inequality $w(v) > \alpha_w(G[N(v)\backslash (IC \cup N(u))])$ is satisfied.  
    \item[$(b)$] Assume that set $C \vartriangleleft A_C$, then $\forall p \in C, w(p) < \alpha_w(G[N(p)\backslash C ])$ is always satisfied. In addition, if there exists a vertex $v \in C$ and a vertex $u \in N^2(v)$ such that $w(v) \geq  \alpha_w(G[N(v) \backslash (C \cup N(u))])$, then $C \cup \{ u \} \vartriangleleft A_C$.
\end{itemize}
\end{theorem}
\proof{Proof} 
We first consider $(a)$ and let set $I_w = V \backslash VC_w$. From the upper bound lemma, these results can be directly obtained: $\forall p \in IC, w(p) \leq w(I_w \cap N[p]) =  w(I_w \cap N(p)) \leq \alpha_w(G[N(p)\backslash IC])$. Also, based on the assumption about $u$ in $(a)$, if $u \in I_w$, then $w(v) \leq w(I_w \cap N[v]) = w(I_w \cap N(v)) \leq \alpha_w(G[N(v)\backslash (IC \cup N(u))])$, which leads to a contradiction.

Similar methods can be used to prove $(b)$. First, $\forall p \in C, \forall I \in A_I, w(p) < w(I \cap N[p]) = w(I \cap N(p)) \leq \alpha_w(G[N(p)\backslash C ])$ can be obtained from the upper bound lemma. Besides, under given conditions about $u$ in $(b)$, if there is an $I^{*} \in A_I$ such that $u \in I^{*}$, a contradiction is deduced from $w(p) < w(I^{*} \cap N[p]) = w(I^{*} \cap N(p)) \leq \alpha_w(G[N(p)\backslash (C \cup N(u)) ])$.
\endproof

The vertex cover extension theorem describes how to expand a set that is~(strongly) sheathed: Attempt to find a vertex that satisfies the condition that after removing its neighbor set, the upper bound lemma is not satisfied in the local structure of the expanded set. If such a vertex is found, it is directly added to the expanded set.

\section{Causal Inference Techniques}
\label{sec:Causal Inference Technique}

In this section, with the help of the upper bound lemma and two extension theorems, we give the CITs used in this paper. Our CITs can be divided into two types: The first type is a strongly reducible state-preserving technique introduced in Section~\ref{subsec:strongly reducible state-preserving technique}, while the second type is a reducible state-preserving technique shown in Section~\ref{subsec:reducible state-preserving technique}.

\subsection{Strongly reducible state-preserving technique}
\label{subsec:strongly reducible state-preserving technique}
The strongly reducible state-preserving technique exploits the assumption that a vertex is strongly reducible, and the assumed state of the vertex can be divided into two cases: The vertex is assumed to be strongly inclusive or is assumed to be strongly sheathed. We first consider the assumption that a vertex is strongly inclusive and give the following definition.
\begin{definition}
\label{de:strongly inclusive assumption}
Let set $S$ be an independent set in the graph. If a vertex $u \in N(S)$ such that $w(u) \geq w(S \cap N(u))$, we call it a child of set $S$. A child $u$ is called an extending child if and only if there exists a unique independent set $S^{*} \subseteq N(u) \backslash N[S]$ such that $ w(u) < w(S \cap N(u)) + w(S^{*})$ and vertex set $S^{*}$ is called a satellite set of set $S$.  
\end{definition}

On the basis of Definition~\ref{de:strongly inclusive assumption}, with the assumption that a vertex is strongly inclusive, the concept of `confined/unconfined vertices’ is given by the following conflict analysis process:
\begin{definition}
\label{de:generalize unconfined}
Let $v$ be a vertex in the graph. Suppose set $S := \{v\} \blacktriangleleft A_I$, repeating $(i)$ until $(ii)$ or $(iii)$ holds:
\begin{itemize}
    \item[$(i)$] As long as set $S$ has an extending child in $N(S)$, set $S$ is extended by including the corresponding satellite set into set $S$.
    \item[$(ii)$] If a child $u$ such that $w(u) \geq w(S \cap N(u)) + \alpha_w(G[N(u) \backslash N[S]])$ could be found, that is, the upper bound lemma is not satisfied in the local structure of set $S$, then halt and vertex $v$ is called an unconfined vertex.
    \item[$(iii)$] If any child is not an extending child, then halt and return set $S_v = S$. In this case, vertex $v$ is called a confined vertex and the set $S_v$ is called the confining set of vertex  $v$.
\end{itemize}
\end{definition}
\begin{figure}[ht]
    	\includegraphics[width=0.70\textwidth]{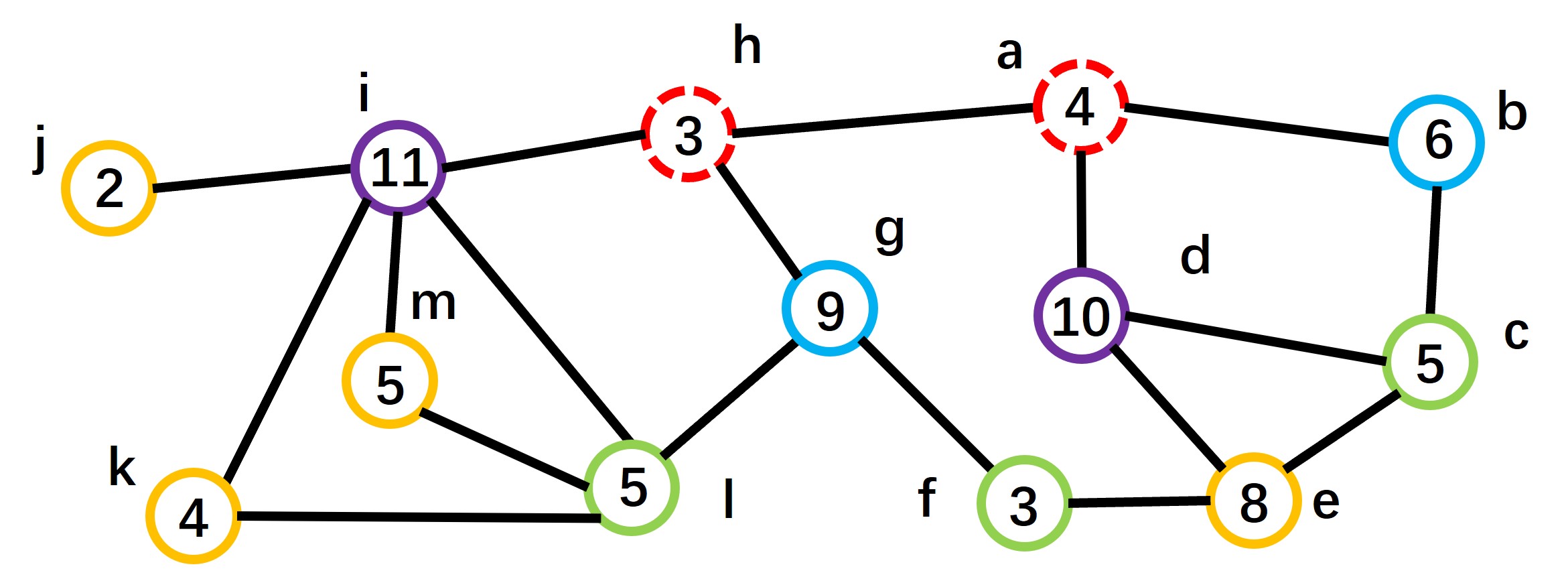}
    	\centering
    	\caption[caption]{Some examples of unconfined vertices, and a MWIS in this graph is $\{ b, d, g, i\}$.  Let set $S := \{ a \}$, from Definition~\ref{de:strongly inclusive assumption}, vertex $b$ is an extending child of set $S$ and set $\{c\}$ is a satellite set of set $S$. Thus, set $S$ can be extended as: $\{a, c\}$. At this time, it can be found that a child $d$ such that $w(d) \geq w(S \cap N(d)) + \alpha_w(G[N(d) \backslash N[S]])$, then halt and conclude that vertex $a$ is an unconfined vertex. Similarly, let set $S := \{ h \}$, then it can be found that vertex $g$ is an extending child of set $S$ and set $\{f,l\}$ is a satellite set of set $S$. So set $S$ can be further expanded as: $\{h, l, f\}$. After that, the child $i$ satisfied: $w(i) \geq w(S \cap N(i)) + \alpha_w(G[N(i) \backslash N[S]])$, hence, vertex $h$ is an unconfined vertex.}
\label{fig:unconfined vertex example}
\end{figure}

Some examples of unconfined vertex are given in Figure~\ref{fig:unconfined vertex example}. By means of the conflict analysis process in Definition~\ref{de:generalize unconfined}, vertices $a$ and $h$ can be found to be unconfined vertices. It is also worth noting that, by the definition of unconfined vertex given in~\cite{xiao2021efficient}, in Figure~\ref{fig:unconfined vertex example}, only vertex $a$ can be found to be an unconfined vertex. The reason for this is that we further generalize the concept of confined/unconfined vertices in this work. Compared with the definition of extending child $u$ in~\cite{xiao2021efficient}, which requires $\vert N(u) \backslash N[S] \vert $ $= 1$ and $w(u) < w(N(u)\backslash N(S))$,  we can consider the more general case where $N(u) \backslash N[S]$ is an independent set rather than a single vertex, helping us find more unconfined vertices.

Next, we will explore the properties of confined/unconfined vertices. By the conflict analysis process in Definition~\ref{de:generalize unconfined} and the independent set extension theorem, set $S$ can be extended under the assumption: set $S := \{v\} \blacktriangleleft A_I$, and set $S \blacktriangleleft A_I$ is always satisfied. If vertex $v$ is a unconfined vertex, then the upper bound lemma is not satisfied in the local structure of set $S$, which contradicts set $S \blacktriangleleft A_I$. Thus, vertex $v$ is sheathed. Otherwise, then there is a state-preserving result, i.e., the corresponding confining set $S_v \blacktriangleleft A_I$ holds. Furthermore, suppose two confined vertices $u$, $v$ and the corresponding confining sets $S_u, S_v$ such that $u \in S_v$ and $v \in S_u$. If $\{v\} \blacktriangleleft A_I$, then obviously $\{u\} \blacktriangleleft A_I$ holds. If not, vertex $v$ is sheathed in graph $G$. Since $v \in S_u$, then vertex $v$ is included in the satellite set of an intermediate state set $S^{\prime}$ of $S_u$, which means that in graph $G[V \backslash \{v\}]$, the upper bound lemma is not satisfied in the local structure of set $S^{\prime}$. Thus, by Definition~\ref{de:generalize unconfined}, vertex $u$ is an unconfined vertex of graph $G[ V \backslash \{v\}]$ and is sheathed in this graph. From these analysis results and the symmetry of the relationship between vertex $v$ and vertex $u$, we can know that vertex set $\{u, v\}$ is a simultaneous set. Therefore, the following properties can be obtained:
\begin{corollary}
\label{pro:Corollary of unconfined vertex and confined vertex}
Let $v$ is a vertex in the graph.
\begin{itemize}
    \item[$(a)$] If vertex $v$ is an unconfined vertex, then it is sheathed and after deleting it from the graph, the weight of the MWIS in the remaining graph remains unchanged.
    \item[$(b)$] Suppose vertex $v$ is a confined vertex, then either it is sheathed or the corresponding confining set $S_v \blacktriangleleft A_I$. Moreover, if a vertex $u \in S_v$ is also a confined vertex with the corresponding confining set $S_u$ and $v \in S_u$, then vertex set $\{u, v\}$ is a simultaneous set.
\end{itemize}
\end{corollary}

From Corollary~\ref{pro:Corollary of unconfined vertex and confined vertex}, it can be known that the conflict analysis process in Definition~\ref{de:generalize unconfined} can be used to find the vertex that is sheathed  or a simultaneous set. These CITs will be used to design reduction rules in Section~\ref{subsec:The Causal Reduce}. In addition, by the property of confined vertex, a fact is obvious: If confined vertex $v$ such that $ \{ v\} \blacktriangleleft A_I$, then the corresponding confining set $S_v \blacktriangleleft A_I$. We will exploit this state-preserving result in the B\&R algorithm to design a branching rule to search for a solution in Section~\ref{subsec:The Causal B&R Solver}. 

Next, we proceed to consider the assumption that a vertex is strongly sheathed. In the MIS problem, the notion of mirror is given by means of such an assumption and is very useful in practice~\cite{akiba2016branch}. We will generalize the notion of mirror to the MWIS problem: For a vertex $v \in V$, a mirror of vertex $v$ is a vertex $u \in N^2(v)$ such that $w(v) \geq \alpha_w(G[N(v) \backslash N(u)])$. 
\begin{remark}
When the weight of all vertices in the graph is $1$, then $\alpha(G[N(v) \backslash N(u)]) = \alpha_w(G[N(v) \backslash N(u)])$ $\leq w(v) = 1$. This means that $N(v) \backslash N(u)$ induces a clique or is an empty set, and this is exactly the definition that vertex $u$ is the mirror of vertex $v$ in the MIS problem.
\end{remark}

To make the concept of mirror more practical, we further generalize it to the case of set, which leads to the following definitions:
\begin{definition}
\label{de:strongly sheathed assumption}
Let set $C$ be a vertex subset in the graph. If a vertex $v \in C$  satisfies the inequality: $w(v) < \alpha_w(G[N(v) \backslash C])$, we call it a father of set $C$. Furthermore, if there exists a vertex $u \in N^{2}(v)$ such that $w(v) \geq \alpha_w(G[N(v) \backslash (C \cup N(u))])$, then the father $v$ is called an extending father of set $C$ and vertex $u$ is called a mirror of vertex $v$. We use $M(v)$ to denote the set of mirrors of vertex $v$.
\end{definition}

By means of Definition~\ref{de:strongly sheathed assumption}, and under the assumption that a vertex is strongly sheathed, the concept of `covered/uncovered vertices' is given by the following conflict analysis process:
\begin{definition}
\label{de:uncovered}
Let $v$ be a vertex in the graph. At the beginning, suppose set $C := \{v\} \vartriangleleft A_C$ and repeating $(i)$ until $(ii)$ or $(iii)$ are met:
\begin{itemize}
    \item[$(i)$] When set $C$ has an extending father, extend set $C$ by including the corresponding set of mirrors to set $C$.
    \item[$(ii)$] If there is a vertex $u \in C$ such that $w(u) \geq \alpha_w(G[N(u) \backslash C])$, in this case, the upper bound lemma is not satisfied, then halt and vertex $v$ is called an uncovered vertex.
    \item[$(iii)$]If set $C$ has no extending father, then halt and return set $C_v = C$. In this case, vertex $v$ is called a covered vertex and vertex set $C_v$ is called the covering set of vertex $v$.
\end{itemize}
\end{definition}
\begin{figure}
     \includegraphics[width=0.75\textwidth]{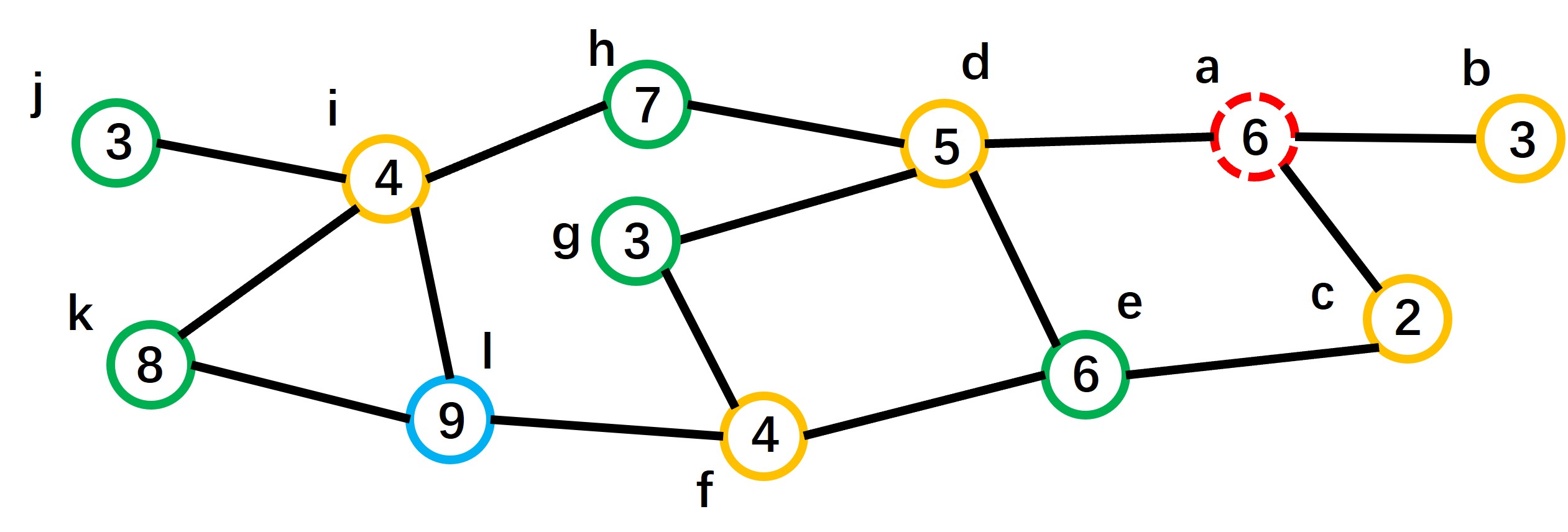}
	\centering
	\caption[caption]{An example of uncovered vertex and a MWIS of this graph is $\{ a, e, g, h,  j, l\}$. Starting with set $C := \{ a \}$, from Definition~\ref{de:strongly sheathed assumption}, it can be seen that vertex $a$ is an extending father of set $C$ and set $\{ e,g,h \}$ is the mirrors set of vertex $a$. Thus, set $C$ can be extended to: $\{ a,e,g,h \}$. Then, vertex $h$ is also an extending father of set $C$ and set $\{ j,k,l \}$ is the mirrors set of vertex $h$. So set $C$ can be further expanded as: $\{ a,e,g,h,j,k,l \}$. At this time, we find that $w(l) \geq \alpha_w(G[N(l) \backslash C])$, then halt and conclude that vertex $a$ is uncovered.}
    \label{fig:uncovered example}
\end{figure}

An example of uncovered vertex is given in Figure~\ref{fig:uncovered example} and we find that vertex $a$ is an uncovered vertex. 
In addition, the properties of uncovered/covered vertices are worth further study. From the vertex cover extension theorem, in the conflict analysis process of Definition~\ref{de:uncovered}, for any extending father $f$ of set $C$, $\forall u \in M(f)$, if set $C \vartriangleleft A_C$, set $C \cup \{ u \} \vartriangleleft A_C$ always holds. Thus, under the assumption set $C := \{v\} \vartriangleleft A_C$, if vertex $v$ is not an uncovered vertex, then a state-preserving result can be obtained: The corresponding covering set $C_v \vartriangleleft A_C$. Otherwise, the upper bound lemma is not satisfied in the local structure of set $C$, which contradicts hypothesis set $C \vartriangleleft A_C$. So vertex $v$ is inclusive. Also, assume that the two covered vertices $u, v$ and the corresponding covering set $C_u, C_v$ satisfy: $v \not \in N(u)$, $u \in C_v$ and $v \in C_u$. If vertex $v$ is inclusive, we first remove $N[v]$ from graph $G$. Since $v \in C_u$, then vertex $v$ is a mirror of an extending father of an intermediate state set $C^{\prime}$ of set $C_u$ and the upper bound lemma cannot be satisfied in graph $G[V \backslash N[v]]$ at this time. Thus, vertex $u$ is an uncovered vertex of graph $G[ V \backslash N[v]]$ and is inclusive in this graph. So there exists a MWIS in graph $G$ containing both vertex $v$ and vertex $u$. Moreover, if $\{v\} \vartriangleleft A_C$, $\{u\} \vartriangleleft A_C$ is clearly satisfied. Thus, from the symmetry of the relationship between vertex $u$ and vertex $v$, it can be known that vertex set $\{u, v\}$ is a simultaneous set. These properties are summarized as follows. 
\begin{corollary}
\label{pro:uncovered and covered}
Let $v$ be a vertex in the graph $G$. 
 \begin{itemize}
     \item[$(a)$] If vertex $v$ is an uncovered vertex, then it is inclusive. After deleting $N[v]$ from the graph, the weight of the MWIS in the remaining graph satisfies: $\alpha_w(G) = \alpha_w(G[V\backslash N[v]]) + w(v)$. 
     \item[$(b)$] If vertex $v$ is a covered vertex. Then, either vertex $v$ is inclusive or the corresponding covering set $C_v \vartriangleleft A_C$. Also, if another covered vertex $u$ with the corresponding covering set $C_u$ satisfies: $v \not \in N(u)$, $u \in C_v$ and $v \in C_u$, then vertex set $\{u, v\}$ is a simultaneous set.
 \end{itemize} 
\end{corollary}

Corollary~\ref{pro:uncovered and covered} gives the following results: The conflict analysis process in Definition~\ref{de:uncovered} can be applied to find the vertex that is inclusive or a simultaneous set. In Section~\ref{subsec:The Causal Reduce}, we will use these CITs to design reduction rules. Besides, by the property of covered vertex in $(b)$ of Corollary~\ref{pro:uncovered and covered}, we can know a state-preserving result: if the covered vertex $v$ such that $\{v\} \vartriangleleft A_C$, then the corresponding covering set $C_v \vartriangleleft A_C$.

\subsection{Reducible state-preserving technique}
\label{subsec:reducible state-preserving technique}
Similar to the first type of CIT, the reducible state-preserving technique utilizes the assumption that a vertex is reducible, that is, assumes that a vertex is inclusive or sheathed. With these assumptions, we can give state-preserving results similar to the first type of CIT. Before that, we give the following definition.
\begin{definition}
\label{de:assumption}
Let sets $IS$ and $IC$ be two vertex subsets in the graph and set $IS$ is an independent set.
\begin{itemize}
    \item[$(a)$] A vertex $u \in N(IS)$ is called an inferred child of set $IS$ if it holds that $w(u) > w(IS \cap N(u))$. Further, if there is only a unique independent set $IS^{*} \subseteq N(u) \backslash N[IS]$ that satisfies the inequality: $ w(u) \leq w(IS \cap N(u)) + w(IS^{*})$, we call the inferred child $u$ an inferred extending child of set $IS$ and vertex set $IS^{*}$ is called an inferred satellite set of set $IS$.
    \item[$(b)$] A vertex $v \in IC$ is called an inferred father of set $IC$ if it holds that $w(v) \leq \alpha_w(G[N(v) \backslash IC])$. An inferred father $v$ is called an inferred extending father of set $IC$ if there exists a vertex $u \in N^{2}(v)$ such that $w(v) > \alpha_w(G[N(v) \backslash (IC \cup N(u))])$ and vertex $u$ is called an inferred mirror of vertex $v$. Also, $IM(v)$ is used to denote its set of inferred mirrors.
\end{itemize}
\end{definition}

By virtue of Definition~\ref{de:assumption} and the assumption that a vertex is inclusive or sheathed, we can directly give the definitions of inferred confining set and inferred covering set accordingly.

\begin{definition}
\label{de:inferred confining set}
Suppose there are no unconfined vertex in the graph. Let $v$ be a vertex in the graph. Beginning with the assumption set $IS := \{v\} 
\subseteq I_w$. 
\begin{itemize}
    \item[$(i)$] Only if set $IS$ has an inferred extending child in $N(IS)$, set $IS$ can be extended by including the corresponding inferred satellite set to set $IS$.
    \item[$(ii)$] The above process halts if set $IS$ has no inferred extending child in $N(IS)$ and return set $IS_v = IS$. We call vertex set $IS_v$ is the inferred confining set of vertex $v$.
\end{itemize}
\end{definition}

\begin{definition}
\label{de:inferred covering set}
We assume that there are no uncovered vertex in graph. Let $v$ be a vertex in the graph. Starting with the assumption set $IC := \{v\} \subseteq VC_w$.
\begin{itemize}
    \item[$(i)$] While set $IC$ has an inferred extending father, extend set $IC$ by including the corresponding set of inferred mirrors to set $IC$.
    \item[$(ii)$] The above process halts if set $IC$ has no inferred extending father and return set $IC_v = IC$. We call vertex set $IC_v$ is the inferred covering set of vertex $v$.
\end{itemize}
\end{definition}
\begin{figure}[h]
     \includegraphics[width=0.75\textwidth]{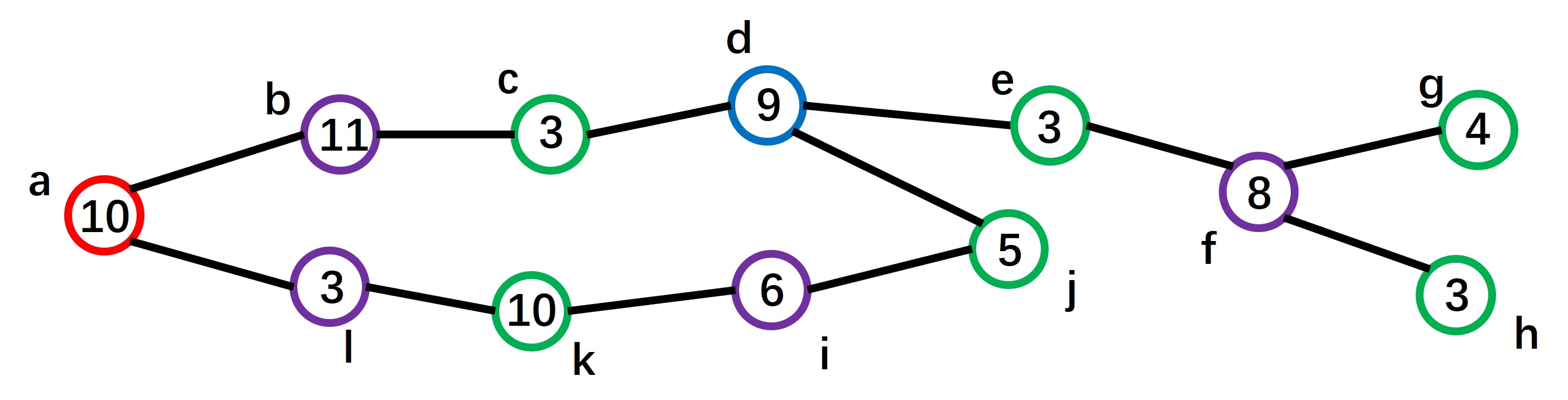}
	\centering
	\caption[caption]{Examples of inferred confining set and inferred covering set. A MWIS for this graph is $\{ a, c, e, g, h, j, k\}$. We first search for the inferred confining set $IS_a$ of vertex $a$. Let set $IS := \{ a \}$, it can be seen from $(a)$ of Definition~\ref{de:assumption} that vertex $b$ is an inferred extending child of set $IS$ and set $\{c\}$ is an inferred satellite set of set $IS$. Thus, set $IS$ can be extended to: $\{a, c\}$. Further, vertex $d$ is also an inferred extending child of set $IS$ and set $\{e, j\}$ is the corresponding inferred satellite set. So set $IS$ can be further extended to: $\{a, c, e, j\}$. At this time, it can be found that both vertex $f$ and vertex $i$ are inferred extending children of set $IS$. Then, both vertex set $\{g, h\}$ and vertex set $\{k\}$ are the corresponding inferred satellite sets. Finally, the inferred confining set of vertex $a$ can be found as: $IS_a = \{a, c, e, j, g, h, k \}$. Furthermore, we continue to search the inferred covering set $IC_d$ of vertex $d$. Let set $IC := \{ d\}$, according to $(b)$ of Definition~\ref{de:assumption}, vertex $d$ is an extending father of set $IC$ and set $\{ b, f, i \}$ is its inferred mirrors set. Then, set $IC$ can be extended as: $\{ b, d , f, i \}$. Next, it can be found that vertex $b$ is an extending father of set $IC$ and set $\{ l \}$ is its inferred mirrors set. Thus, set $IC$ can be further extended as: $\{ b, d , f, i, l \}$. Finally, the inferred covering set of vertex $d$ can be found as: $IC_d = \{ b, d , f, i, l \}$.}
    \label{fig:Inferred graph example}.
\end{figure}

Examples of inferred confining set and inferred covering set are given in Figure~\ref{fig:Inferred graph example}. By the process in Definition~\ref{de:inferred confining set}, we can find the inferred confining set $IS_a = \{a, c, e, j, g, h, k \}$ of vertex $a$. Similarly, according to the process in Definition~\ref{de:inferred covering set}, we can find the inferred covering set $IC_d = \{ b, d , f, i, l \}$ of vertex $d$. Moreover, from the independent set extension theorem and the vertex cover extension theorem, we can directly obtain the following Corollary:
\begin{corollary}
\label{pro:inferred}
Let $v$ be a vertex in the graph. 
\begin{itemize}
    \item[$(a)$] If $\{ v \} \subseteq I_w$, then the corresponding inferred confining set $IS_v \subseteq I_w$.
    \item[$(b)$] Suppose $\{v\} \subseteq VC_w$, then the corresponding inferred covering set $IC_v \subseteq VC_w$.
\end{itemize}
\end{corollary}

From (a) of Corollary~\ref{pro:inferred}, under the premise $\{ v \} \subseteq I_w$, the state-preserving result can be obtained: $IS_v \subseteq I_w$. We will integrate this result into the local search process of heuristic algorithm in Section~\ref{subsec:The Causal Search}. In addition, (b) of Corollary~\ref{pro:inferred} also gives a similar state-preserving result result: If $\{ v \}\subseteq VC_w$, then the corresponding inferred covering set $IC_v \subseteq VC_w$. This result can be used to design a branching rule to search for a solution in Section~\ref{subsec:The Causal B&R Solver}.

Furthermore, during the branching process of the B\&R algorithm, it is assumed that a vertex $v$ is selected for branching. Inspired by the successful application of packing constraints in the MIS problem, we extend them to the MWIS problem and propose the concept of ``weight packing constraint''. 

When assuming that vertex $v$ is inclusive, if $ \exists u \in N(v)$ such that $w(u) \geq w(v)$, let $N^{+}(u) = N(u)\backslash N[v]$. To avoid obtaining another MWIS by adding vertex $u$ to the independent set and removing vertices in $N(u)$ from the independent set, by the upper bound lemma, the following state-preserving result needs to be guaranteed to hold: 
$$ w(v)  + \sum\limits_{z \in N^{+}(u)} w(z)(1 - x_z) > w(u).$$
The $0$-$1$ integer variable $x_z$ is used to indicate whether vertex $z \in N^{+}(u)$ is in the independent set, and $x_z = 0$ means it is in the independent set, otherwise it is not. Thus, a weight packing constraint can be created as shown below:
\begin{equation}
\label{eq:packing constraint in MWIS}
 \sum\limits_{z \in N^{+}(u)} w(z) x_z < \sum\limits_{z \in N^{+}(u)} w(z) - ( w(u) - w(v)).
\end{equation}

When assuming that vertex $v$ is sheathed, to avoid that a MWIS containing it can be found by modifying its state, by means of the upper bound lemma, the following state-preserving result needs to be satisfied:
$$\sum\limits_{z \in N(v)} w(z) (1 - x_z) > w(v).$$ 
So a weight packing constraint can also be created as follows:
\begin{equation}
\label{eq:packing constraint not in MWIS}
\sum_{z \in N(v)} w(z) x_z < \sum_{z \in N(v)} w(z) - w(v).
\end{equation}

These constraints will be kept and managed while the algorithm is searching for a solution, and we only need to search all branches satisfying these constraints, since no better solution exists in the remaining branches, thus narrowing the search space. Let $\sum\limits_{z \in S} w(z) x_z < k$ be a weight packing constraint such that set $S$ is non-empty. When a vertex $z$ is found to be inclusive, for each constraint that includes variable $x_z$, we delete the variable on the left side of the constraint and keep the right side of the constraint unchanged. When a vertex $z$ is inferred to be sheathed, for each constraint that contains variable $x_z$, we delete the variable on the left side of the constraint and decrease the weight of vertex $z$ on the right side of the constraint. In the process of keeping and managing these constraints, some properties of causal inference are mined, which can be divided into the following three cases.
\begin{itemize}
    \item[$(a)$] When there is a constraint whose right-hand term $k$ is less than or equal to $0$, then we can directly prune subsequent searches from the current branch vertex.
    \item[$(b)$] When there is a constraint whose right-hand term $k$ is less than or equal to the weight of any vertex in set $S$, if this set is not an independent set, we can prune subsequent searches from the current branch vertex. If not, the vertices in set $S$ will be included in the independent set. 
    
    In addition, some new weight packing constraints can also be introduced. Suppose there is a vertex $p \in N(S)$ such that $w(p) \geq w(N(p) \cap S)$, let $N^{+}(p) = N(p) \backslash N[S]$, by the upper bound lemma, the following state-preserving result needs to be guaranteed:
    $$w(N(u) \cap S) + \sum\limits_{z \in N^{+}(u)} w(z) (1 - x_z) > w(u).$$
    Therefore, we can introduce the following weight packing constraint:
    \begin{equation}
        \label{eq:packing set in MWIS}
        \sum\limits_{z \in N^{+}(p)} w(z) x_z < \sum\limits_{z \in N^{+}(p)} w(z) - ( w(p) - w(N(p) \cap S)).
    \end{equation}
    \item[$(c)$] When there is a constraint whose right-hand term $k > 0$ and there is vertex $u \in N(S)$ such that $\sum\limits_{z\in N(u) \cap S} w(z)$ $ \geq k$, it can be inferred that vertex $u$ is sheathed to ensure that this constraint holds. In addition, in order to ensure that the current state-preserving result is valid, similar to constraint~\eqref{eq:packing constraint not in MWIS}, the following constraint needs to be introduced:
    \begin{equation}
    \label{eq:packing set constraint not in MWIS}
    \sum_{z \in N(u)} w(z) x_z < \sum_{z \in N(u)} w(z) - w(u).
    \end{equation}
\end{itemize}
The above properties of causal inference provide new pruning search techniques for the B\&R algorithm and can simplify the graph. We will integrate these techniques into B\&R algorithm in Section~\ref{subsec:The Causal B&R Solver}.

\section{Integrate CITs into Existing Algorithmic Frameworks}
\label{sec:Integrate CITs into Existing Algorithmic Frameworks}
We next describe in detail how CITs in Section~\ref{sec:Causal Inference Technique} are integrated into the existing algorithmic frameworks. Section~\ref{subsec:The Causal Reduce} introduces how to apply the first type of CIT to the reduction algorithm. Further, integrating the resulting reduction algorithm and the state-preserving results of two types of CITs into B\&R algorithm will be presented in Section~\ref{subsec:The Causal B&R Solver}, and Section~\ref{subsec:The Causal Search} will introduce the application of the state-preserving results of the second type of CIT to the local search process of heuristic algorithm.
\subsection{The Causal Reduce}
\label{subsec:The Causal Reduce}
We first introduce how to design reduction rules with the first type of CIT and how to integrate them into the existing reduction algorithm. From the property of unconfined vertex in Corollary~\ref{pro:Corollary of unconfined vertex and confined vertex} and the property of uncovered vertex in Corollary~\ref{pro:uncovered and covered}, the following reduction rules that can directly determine whether a vertex is reducible are given first:
\begin{itemize}
    \item \textbf{Rule I:} Check whether a vertex $v$ is unconfined or confined by the procedure in Definition~\ref{de:generalize unconfined}, and if it is unconfined, remove vertex $v$ directly from the graph.
    \item \textbf{Rule II:} Use the procedure in Definition~\ref{de:uncovered} to check whether a vertex $v$ is covered or uncovered, and if it is uncovered, include vertex $v$ into the independent set and remove $N[v]$ from the graph.
\end{itemize}

Before further introducing how to utilize the first type of CIT to design reduction rules, we first give an important property about simultaneous set mentioned in~\cite{xiao2021efficient}:  A simultaneous set $S$ can be contracted by removing all vertices in set $S$ from the graph and introducing a vertex $v^{*}$ such that it is adjacent to all vertices in $N(S)$ with weight $w(v^{*}) = w(S)$, while the weight of the MWIS in the remaining graph remain unchanged. 

Next, we will design reduction rules on simultaneous set through the first type of CIT, and give the following definitions by the results of the simultaneous set given in $(b)$ of Corollary~\ref{pro:Corollary of unconfined vertex and confined vertex} and $(b)$ of Corollary~\ref{pro:uncovered and covered}.
\begin{definition}
\label{de:simultaneous set}
Let $u, v$ be two vertices in the graph.
\begin{itemize}
    \item[$(a)$] Suppose  vertices $u$ and $v$ be two confined vertices with confining set $S_u$ and $S_v$. If $u \in S_v$ and $v \in S_u$, then set $\{ u, v \}$ is called a confining simultaneous set.
    \item[$(b)$] Assume that vertices $u$ and $v$ be two covered vertices with covering set $C_u$ and $C_v$. Set $\{u, v\}$ is called a covering simultaneous set if $u \in C_v$ and $v \in C_u$.
\end{itemize}
\end{definition}

From Definition~\ref{de:simultaneous set}, we have the following rules:
\begin{itemize}
    \item \textbf{Rule III:} If there are two confined vertices that constitute a confining simultaneous set, then merge them.
    \item \textbf{Rule IV:} Merge two covered vertices $u$ and $v$ if they form a covering simultaneous set.
\end{itemize}

Next, we will describe how to integrate our reduction rules into an existing reduction algorithm---\textbf{Reduce} proposed by~\cite{xiao2021efficient}. \textbf{Reduce} consists of seven steps. The reduction rules used in these steps exploit the sufficient conditions that a vertex is reducible. It executes these steps incrementally, which means that the next step is only executed when all previous steps are no longer applicable. Thus, if the graph is changed, it will go back to the first step. Notably, our reduction rules \textbf{I} and \textbf{III} are further generalization of the reduction rules used in step~$5$ of \textbf{Reduce}. So, we can combine our reduction rules \textbf{I} and \textbf{III} into one step to replace step~$5$ in \textbf{Reduce} and label this step as \textbf{Remove Unconfined \& Contract Confining}. Similarly, we can also integrate our reduction rules \textbf{II} and \textbf{IV} into another new step in the reduction algorithm, called \textbf{Remove Uncovered \& Contract Covering}.
\begin{itemize}
    \item \textbf{Remove Unconfined \& Contract Confining:} Check whether a vertex is unconfined or confined. If it is confined, apply \textbf{Rule I} to remove it; If not, use \textbf{Rule III} to contract the corresponding confining simultaneous set when it can be found.
    \item \textbf{Remove Uncovered \& Contract Covering:} If a vertex is checked to be uncovered, use \textbf{Rule II} to reduce it. Otherwise, if the corresponding covering simultaneous set can be found, use \textbf{Rule IV} to merge it.
\end{itemize}

\begin{figure}[h]
     \includegraphics[width=0.90\textwidth]{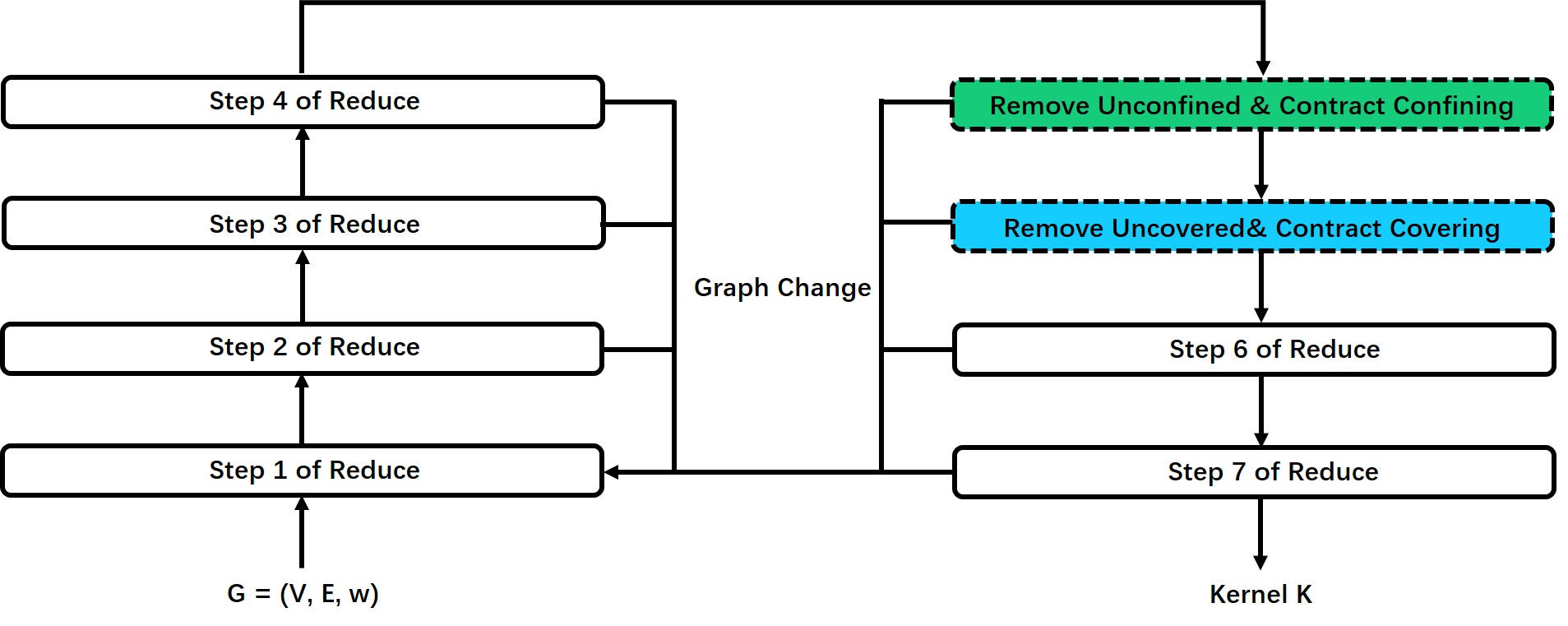}
	\centering
	\caption[caption]{Casual Reduce: Given an input graph $G$, each step of the algorithm is executed sequentially and the graph changes, immediately go back to the first step. When all steps are completed and the graph no longer changes, return to the remaining graph \textbf{kernel}.}
    \label{fig:Causal Reduce}
\end{figure}

Thus, a new reduction algorithm called \textbf{Causal Reduce} can be obtained by using \textbf{Remove Unconfined \& Contract Confining} to replace step $5$ of \textbf{Reduce} and adding \textbf{Remove Uncovered \& Contract Covering} between \textbf{Remove Unconfined \& Contract Confining} and step $6$ of \textbf{Reduce}, which is shown in Figure~\ref{fig:Causal Reduce}. We will use \textbf{Causal Reduce}$(G)=(K, c)$ to represent the processing of this algorithm on a given input graph $G$. The processing result of this algorithm consists of two parts: One is the remaining graph called kernel $K$ and the other is the weight of the vertex set contained in the MWIS obtained by inference. It’s worth noting that the reduction algorithm \textbf{Causal Reduce} may not resolve all instances directly, but it can be used as a preprocessing for heuristic and exact algorithm.
\subsection{The Causal B\&R Solver}
\label{subsec:The Causal B&R Solver}
Before introducing how to integrate our CITs into B\&R algorithm, we briefly introduce the state-of-the-art exact algorithm \textbf{Solve} proposed by~\cite{xiao2021efficient}. \textbf{Solve} is based on the idea of B\&R algorithm, which first apply reduction algorithm \textbf{Reduce} to reduce the instance. Then, apply branching rule by virtue of the property of the confining set and perform reduction algorithm \textbf{Reduce} in every branch of the search tree to find a solution. During the searching, it uses a standard technique based on finding upper and lower bounds to prune the search tree and take the best solution weight $W_b$ currently found in the algorithm as the lower bound. Initially, let $W_b$ be the weight of the solution obtained by heuristic algorithm on the kernel $K$, and update $W_b$ once a better solution is obtained in the algorithm. The heuristic algorithm, denoted by \textbf{Greedy}($G$), is a greedy algorithm that iteratively selects a vertex in order of some measure and removes its closed neighbor set from the graph. In each searching branch, it uses a heuristic method to find an upper bound $W_{ub}$ of the optimal solution weight of the current graph, which is based on weight clique covers and is denoted by \textbf{UpperBound}($G$).  If the current best solution weight $W_b$ is not smaller than $W_{ub}$, then there is no better solution in this searching branch and it can be discarded directly. 
\begin{algorithm}[ht]
\caption{The \textbf{Causal B\&R Solver}($G$)}
\label{alg:The exact algorithm}
\begin{algorithmic}[1]
\scriptsize
\REQUIRE A vertex weight graph $G = (V, E,w)$;
\ENSURE The weight of a MWIS of $G$.
\STATE Initialization of global variable $W_b$: $W_b \leftarrow 0$;
\IF{weight packing constraints have been created}
\WHILE{True}
\STATE  $(K, c) \leftarrow\textbf{Causal Reduce}(G)$;
\STATE   \textbf{check constraints}();
\IF{existence constraints are not satisfied}
\STATE \textbf{return} $W_b$;
\ELSIF{graph is simplified}
\STATE continue;
\ELSE
\STATE break;
\ENDIF
\ENDWHILE
\ELSE 
\STATE $(K, c) \leftarrow\textbf{Causal Reduce}(G)$;
\ENDIF
\STATE  $W_b \leftarrow \max\{W_b, c + \textbf{Greedy}(K)\}$;
\IF{$c + \textbf{UpperBound}(G) \leq W_b$}
\STATE \textbf{return} $W_b$;
\ENDIF
\STATE Pick up a vertex $v$ of maximum degree and compute the confining set $S_v$ and the inferred covering set $IC_v$;
\STATE create weight packing constraint~\eqref{eq:packing constraint in MWIS} and $W_b \leftarrow \max\{W_b, c + w(S_v) + \textbf{Causal B\&R Solver}(K - N[S_v])\}$;
\STATE create weight packing constraint~\eqref{eq:packing constraint not in MWIS} and $W_b \leftarrow \max\{W_b, c + \textbf{Causal B\&R Solver}(K - IC_v ) \}$;
\STATE \textbf{return} $W_b$;
\end{algorithmic}
\end{algorithm}

Our CITs will be integrated into two parts of \textbf{Solve}, resulting in a new exact algorithm called \textbf{Causal B\&R Solver}. The first part is that we will use our reduction algorithm \textbf{Causal Reduce} to reduce the instance to get the kernel $K$, and perform the reduction algorithm on each branch of the search tree. The second part is that we will make use of the state-preserving results of two types of CITs during the branching process. Similar to the idea of \textbf{Solve} in \cite{xiao2021efficient}, using property of confining set to the branching process, when choosing a vertex with the maximum degree to branch, the state-preserving results of $(b)$ of Corollary~\ref{pro:inferred} and $(b)$ of Corollary~\ref{pro:Corollary of unconfined vertex and confined vertex} will be used in this part. This means that during branching, we either remove the inferred covering set of the branching vertex from the graph or include the confining set of the branching vertex into the independent set. Furthermore, we will create weight packing constraint~\eqref{eq:packing constraint not in MWIS} while removing the inferred covering set of branching vertex. Similarly, we will also create weight packing constraint~\eqref{eq:packing constraint in MWIS} when including the confining set of branching vertex into the independent set. We will keep and manage these weight packing constraints when searching for solutions in each branch of the search tree. Specifically, another step called \textbf{check constraints} is added after the last step of \textbf{Casual Reduce}. In this step, for each weight packing constraint, we will check whether the constraint holds and whether the graph can be simplified by the causal inference properties of that constraint. If any constraint is violated, the searching branch will be skipped. If the graph can be simplified, \textbf{Causal Reduce} will continue to execute after reducing the graph. If none of the above conditions are met, the subsequent process will be performed. The main steps of \textbf{Causal B\&R Solver} are listed in Algorithm~\ref{alg:The exact algorithm}.
\subsection{The Causal Search}
\label{subsec:The Causal Search}
After taking our reduction algorithm \textbf{Causal Reduce} as a preprocessing, we apply the state-preserving result of second type of CIT to the local search process of heuristic algorithm DynWVC2~\cite{cai2018improving} to solve the complementary problem of the MWIS problem---the MWVC problem, which leads to a new algorithm called \textbf{Causal Search}.
\begin{algorithm}[h]
\caption{The basic framework of DynWVC2 algorithm.}
\label{alg:The basic framework of the DynWVC2 algorithm for solving MWVC}
\begin{algorithmic}[1]
\scriptsize
\REQUIRE A vertex weight graph $G = (V, E,w)$, the cutoff time of the running $T$;
\ENSURE A vertex cover of $G$.
\STATE  $VC \leftarrow$ \textbf{Construct}();
\STATE  $VC^{*} \leftarrow VC$;
\WHILE{$elapsed\_time \textless T$}
\STATE $R \leftarrow$ \textbf{RemoveVertices}(VC)
\WHILE{some edge is uncovered by $VC$}
\STATE choose a vertex $v$ from $N(R)$;
\STATE  $VC \leftarrow VC \cup \{ v \}$;
\ENDWHILE
\STATE remove redundant vertices from $VC$;
\IF{$w(VC) \textless w(VC^{*})$}
\STATE $VC^{*} \leftarrow VC$
\ENDIF
\ENDWHILE
\end{algorithmic}
\end{algorithm}

The DynWVC2 algorithm proposed by \cite{cai2018improving}, is the state-of-the-art heuristic algorithm for solving MWVC problem. The basic framework of this algorithm is shown in Algorithm~\ref{alg:The basic framework of the DynWVC2 algorithm for solving MWVC}. The local search process of this algorithm mainly consists of a removing phase and an adding phase, and the specific process can be found in \cite{cai2018improving}.  

Our CITs will be considered in the removing phase of the algorithm --- \textbf{RemoveVertices} function. In this function, there are two scoring functions $loss$ and $valid\_score$ used to select the vertices to remove from the vertex cover $VC$. The specific definition of these two scoring functions can be seen in~\cite{cai2018improving}. The $loss$ and $valid\_score$ functions have fundamentally different effects on the behavior of the algorithm. Vertex selection using $loss$ function is an ``exploratory” selection; in other words, it is quite possible that such a chosen vertex is good for the quality of the solution, but this cannot be determined. Different from ``exploratory” vertex selection, $valid\_score$ is a ``deterministic” selection, that is, we can determine whether removing a vertex will have a positive impact on the quality of the solution. For example, if a vertex has a negative $valid\_score$ value, this means that after removing this vertex and adding its adjacent uncovering vertices, a vertex cover with lower weight than the current vertex cover can be obtained~\cite{cai2018improving}.

In removing phase, the vertex with the minimum loss is removed from vertex cover $VC$ first, and then the second removed vertex is selected by a dynamic vertex selection strategy. The details of dynamic vertex selection strategy can be learned in~\cite{cai2018improving}. After removing the two vertices, if the total degree of the removed vertices does not reach a predetermined value~(which is set to $2$ times average degree of the graph), another vertex to be selected with the BMS strategy~\cite{cai2015balance}, which samples $t$~($t=50$) vertices from vertex cover $VC$ and chooses the one with the minimum $loss$, will be removed to expand the search region. In this way, it solves the problem that when removing two vertices the resulting search area is too small and limits the ability of the adding phase to find better local optima. If the search area obtained by removing two vertices is large enough, in order to balance the search time and search quality, the third vertex will not be selected for removing.

\begin{figure}[ht]
     \includegraphics[width=0.95\textwidth]{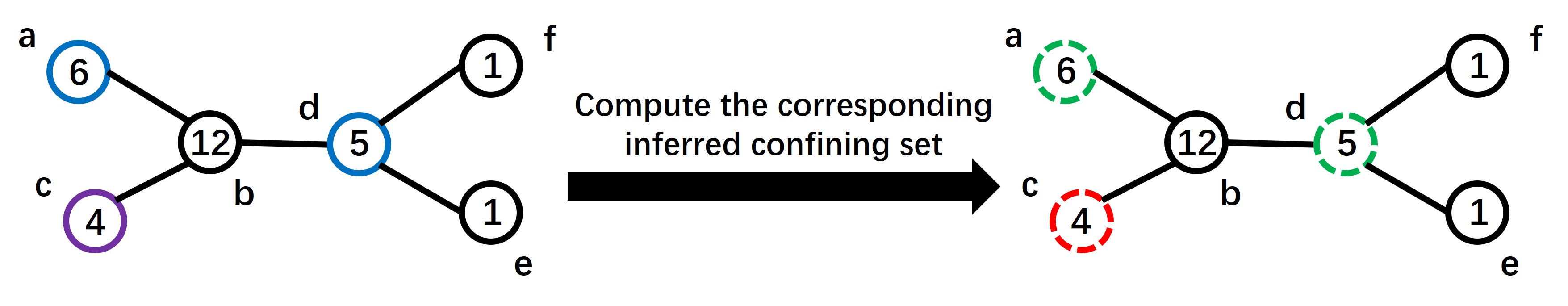}
	\centering
	\caption[caption]{Example of our CITs applied to the vertex removing process: When we utilize $valid\_score$ to select the removed vertex $c$ from the vertex cover $VC = \{ a, c, d\}$, we can compute the corresponding inferred confining set $IS_c = \{a, c, d \}$ of vertex $c$ and remove set $IS_c$ from vertex cover $VC$. }
    \label{fig:eneral local search and local search combined with local causal inference information}
\end{figure}

The state-preserving result of second type of CIT will be applied to the dynamic vertex selection strategy for selecting the second vertex to be removed. The dynamic vertex selection strategy consists of a primary vertex scoring function $valid\_score$ and a secondary scoring function $loss$. When the removed vertex $v$ is selected by $valid\_score$ function, it can be seen from the nature of the $valid\_score$ function: There is a high probability that there exists a MWIS $I$ containing it. If the vertex $v$ is indeed included in $I$, by $(a)$ of Corollary~\ref{pro:inferred}, the corresponding inferred confining set $IS_v$ also contained in $I$. Inspired by this result, when selecting the second removed vertex $v$ by scoring function $valid\_score$, we will remove the vertices in the inferred confining set $IS_v$ from the vertex cover $VC$.
In this way, the search region can be expanded and the number of times to continue to use the third removed vertex to expand the search area is reduced, which means that the ability of local search to find better local optima is improved. An example of our CITs applied to the vertex removing process is presented in Figure~\ref{fig:eneral local search and local search combined with local causal inference information}. 

In addition, it can be seen from  the calculation process of Definition~\ref{de:inferred confining set} about the inferred confining set: the computational complexity of $IS_v$ for each vertex $v$ is $O(\vert N(IS_v) \vert \vert IS_v \vert)$. This means that in the actual application process, since the size of the generally obtained inferred confining set is relatively small, its computational cost is very small. Thus, our CITs is helpful for improving the performance of local search process.

\section{Experiments}
\label{sec:Experiments}
We will conduct four experiments to verify the effect of integrating our CITs into current algorithmic frameworks. The first experiment is used to analyze the impact of our CITs for the reduction algorithm. The examination of the performance gain of our CITs in the B\&R algorithm is shown in the second experiment. The third experiment is used to test the ability of our \textbf{Causal Reduce} as a preprocessing to improve the performance of the heuristic algorithm. The last experiment is conducted to verify the effect of adding our CITs to the local search process of heuristic algorithm.

\noindent\textbf{Experiment environment Setup.}
All of our algorithms are implemented in C++, and compiled by g++ with `-O3’ option. All experiments are run on a platform with $128$G RAM and one Intel(R) Xeon(R) Gold $5117$ CPU @ $2.00$GHz.

\noindent\textbf{Compared Algorithms.}
In previous studies, most of them only use some simple rules as preprocessing to reduce problem instances, and do not pay attention to the performance of preprocessing. Two recent papers~\cite{lamm2019exactly,xiao2021efficient} have studied in depth the reduction rules for the MWIS and analyzed their performance. Since the algorithm \textbf{Reduce} in~\cite{xiao2021efficient} outperforms the algorithm in~\cite{lamm2019exactly} and our \textbf{Causal Reduce} is obtained by integrating our CITs into \textbf{Reduce}, in this paper, we only use it as a baseline to analyze the impact of our CITs for the reduction algorithm. Additionally, in order to fully understand the role of different CITs on the reduction algorithm, we control the application of CITs in \textbf{Reduce} and conduct comparative experiments. Similar to the Causal Reduce shown in Figure~\ref{fig:Causal Reduce}, we use \textbf{Re-Confin} to represent the algorithm obtained after replacing the step~$5$ of \textbf{Reduce} with \textbf{Remove Unconfined \& Contract Confining} and \textbf{Re-Cover} to denote the algorithm obtained by adding \textbf{Remove Uncovered \& Contract Covering} between step~$5$ and $6$ of \textbf{Reduce}.

On the basis of the reduction algorithm \textbf{Reduce}, the authors of~\cite{xiao2021efficient} also developed a fast exact algorithm~\textbf{Solve}, which is the state-of-art exact algorithm in previous work, and our \textbf{Causal B\&R Solver} is obtained by applying our CITs into it, so it will be used as a baseline to verify the performance improvement of our CITs for the B\&R algorithm. Furthermore, we use \textbf{Solve-CR} to identify the algorithm obtained by replacing \textbf{Reduce} with \textbf{Causal Reduce} in \textbf{Solve}, \textbf{Solve-CR-IC} refers to the algorithm obtained by further simplifying the branch by using the inferred covering set of the branching vertex in the branching process on the basis of \textbf{Solve-CR}, and \textbf{Solve-Packing} to represent the algorithm obtained by applying our weight packing constraints to the branching process of \textbf{Solve}.  We will conduct comparative experiments on these algorithms to clarify the impact of different CITs on the B\&R algorithm.

Two state-of-the-art heuristic algorithms FastWVC~(\textbf{Fast})~\cite{cai2019towards} and DynWVC2~(\textbf{Dyn})~\cite{cai2018improving} will be used to verify that our \textbf{Causal Reduce} as preprocessing improves the performance of the heuristic algorithm. We will use \textbf{Causal Re + Fast} and \textbf{Causal Re + Dyn} to denote applying our \textbf{Causal Reduce} as preprocessing before executing FastWVC and DynWVC2. In addition, to further verify the superiority of our \textbf{Causal Reduce} as preprocessing for improving the performance of the heuristic algorithm, we also conduct comparative experiments using \textbf{Reduce} as a preprocessing of the heuristic algorithm. Likewise, we use \textbf{Re + Fast} and \textbf{Re + Dyn} to indicate the application of the previous reduction algorithm \textbf{Reduce} before FastWVC and DynWVC2 are executed.  Moreover, our \textbf{Causal Search} is obtained by integrating CITs into the local search process of DynWVC2. Therefore, we can verify the effect of this operation by comparing DynWVC2 with \textbf{Causal Search}. 

\noindent\textbf{Instances.} We evaluate all algorithms on six real graphs which are most representative and most difficult graphs from different domains. These graphs are downloaded from Network Data Repository~\cite{rossi2015network}. All of them have $100$ thousands to millions of vertices, and dozens of millions of edges. These instances become popular in recent works for the MWIS problem. Statistics of these graphs are shown in Table~\ref{tab:Statistics of real graphs}. In our experiment, the weight of each vertex in the graph will have two random allocation mechanisms~\footnote{All datasets obtained through these two random assignment mechanisms can be found at \url{http://lcs.ios.ac.cn/~caisw/graphs.html}.}, which are commonly used in previous work~\cite{lamm2019exactly,xiao2021efficient,cai2018improving,cai2019towards}. The first allocation mechanism is that the weight of each vertex in the graph is obtained from $[1,200]$ uniformly at random, we will number the six datasets with $1-6$. The second allocation mechanism is that the weight of each vertex in the graph follows a random uniform distribution of $[20,100]$, and $7-12$ will be used to number the six datasets. 
\begin{table}[h]
\centering
\resizebox{\textwidth}{10mm}{
\begin{tabular}{|c|c|c|c|c|c|c|c|c|c|c|c|c|c|c|c|c|}
\hline 
                  &	  \textbf{inf-road-usa}  &  \textbf{soc-livejournal} & \textbf{sc-ldoor}  & \textbf{tech-as-skitter}&  \textbf{sc-msdoor}  & \textbf{inf-roadNet-CA}   \\ \hline 
\textbf{Vertices} & $23947347$                &  $4033137$                &  $952203$          &      $1694616$          & $415863$             &  $1957027$                  \\ \hline
\textbf{Edges}    & $28854312$                &  $27933062$               &  $20770807$        &     $11094209$          & $9378650$            &  $2760388$            \\ \hline
\textbf{NO.}&  $1$, $7$                 &  $2$, $8$                 &  $3$, $9$          &  $4$, $10$              &  $5$, $11$           &  $6$, $12$ \\ \hline
\end{tabular}
}
\caption{All graphs are sorted in descending order regarding the number of edges. In the row headed by ``NO.", each number is used to represent the corresponding number of the dataset generated by the graph according to the corresponding vertex weight allocation mechanism.}
\label{tab:Statistics of real graphs}
\end{table}
\subsection{Impact of CITs on Reduction Algorithm}
\label{subsec:Impact of CITs on Reduction Algorithm}
We first analyze the impact of our CITs for the reduction algorithm and evaluate the performance of all reduction algorithms by measuring the running time, the size of the remaining graphs~(kernel size), and the ratio of the kernel size to the number of vertices in the original graph~(We simply refer to it here as the ratio for convenience.). 
\begin{table}[h]
\centering
\resizebox{\textwidth}{25mm}{
\begin{tabular}{|c|c|c|c|c|c|c|c|c|c|c|c|c|c|c|c|c|}
\hline 
              &                 &  \multicolumn{3}{c|}{\textbf{Reduce}}                  &     \multicolumn{3}{c|}{\textbf{Re-Confin}}           & \multicolumn{3}{c|}{\textbf{Re-Cover}}        &  \multicolumn{3}{c|}{\textbf{Causal Reduce}} \\ \hline
\textbf{NO.}  & $\vert V \vert$ & \textbf{Time(S)}    & \textbf{Kernel Size} & \textbf{Ratio(\%)}& \textbf{Time(S)} & \textbf{Kernel Size} & \textbf{Ratio(\%)} & \textbf{Time(S)}  & \textbf{Kernel Size} &\textbf{Ratio(\%)}  & \textbf{Time(S)}    & \textbf{Kernel Size}  & \textbf{Ratio(\%)} \\ \hline
$\mathbf{1}$  & $23947347$      &   $72.00$           & $431891$      & $1.80$            & $71.66$          & $428137$      & $1.79$             & $74.03$           & $431888$     & $1.80$               & $\underline{55.88}$ & $\bf275082$    & $1.15$             \\ \hline
$\mathbf{2}$  & $4033137$       & $\underline{12.93}$ & $7273$        & $0.18$            & $16.57$          & $5960$        & $0.15$             & $16.20$           & $7261$       & $0.18$               & $17.92$             & $\bf3281$      & $0.08$             \\ \hline
$\mathbf{3}$  & $952203$        & $\underline{3.08}$  & $6492$        & $0.68$            & $3.89$           & $2780$        & $0.29$             & $3.89$            & $6447$       & $0.68$               & $3.84$              & $\bf1682$      & $0.18$             \\ \hline
$\mathbf{4}$  & $1694616$       & $\underline{1.78}$  & $5904$        & $0.35$            & $2.32$           & $5613$        & $0.33$             & $2.26$            & $5909$       & $0.35$               & $5.03$              & $\bf3974$      & $0.23$             \\ \hline
$\mathbf{5}$  & $415863$        & $\underline{1.66}$  & $6570$        & $1.58$            & $1.84$           & $3166$        & $0.76$             & $2.04$            & $6570$       & $1.58$               & $1.86$              & $\bf2162$      & $0.52$             \\ \hline
$\mathbf{6}$  & $1957027$       & $76.14$             & $305470$      & $15.61$           & $63.65$          & $300135$      & $15.34$            & $76.57$           & $305470$     & $15.61$              & $\underline{41.12}$ & $\bf202885$    & $10.37$            \\ \hline
$\mathbf{7}$  & $23947347$      & $117.54$            & $437993$      & $1.83$            & $94.52$          & $434606$      & $1.81$             & $113.07$          & $438004$     & $1.83$               & $\underline{58.02}$ & $\bf235243$    & $0.98$             \\ \hline
$\mathbf{8}$  & $4033137$       & $\underline{11.68}$ & $7620$        & $0.19$            & $15.90$          & $6371$        & $0.16$             & $15.65$           & $7576$       & $0.19$               & $23.38$             & $\bf3806$      & $0.09$             \\ \hline
$\mathbf{9}$  & $952203$        & $8.05$              & $29116$       & $3.06$            & $5.08$           & $7906$        & $0.83$             & $8.88$            & $29116$      & $3.06$               & $\underline{5.06}$  & $\bf4628$      & $0.49$             \\ \hline
$\mathbf{10}$ & $1694616$       & $\underline{1.94}$  & $6999$        & $0.41$            & $2.55$           & $6623$        & $0.39$             & $2.53$            & $6991$       & $0.41$               & $5.19$              & $\bf4429$      & $0.26$             \\ \hline
$\mathbf{11}$ & $415863$        & $4.09$              & $24736$       & $5.95$            & $2.75$           & $10098$       & $2.43$             & $4.24$            & $24736$      & $5.95$               & $\underline{2.66}$  & $\bf7081$      & $1.70$             \\ \hline
$\mathbf{12}$ & $1957027$       & $19.01$             & $131498$      & $6.72$            & $17.00$          & $128637$      & $6.57$            &  $19.54$           & $131498$     & $6.72$               & $\underline{8.64}$  & $\bf65124$     & $3.33$             \\ \hline
\end{tabular}
}
\caption{Impact of CITs for the the reduction algorithm. The \textbf{bold} and \underline{underlined} numbers are the minimum kernel size and shortest running time, respectively.}
\label{tab:Comparison of reduction algorithms performance}
\end{table}   

The experimental results of all algorithms are output in Table~\ref{tab:Comparison of reduction algorithms performance}. We can know that all reduction algorithms can significantly simplify the graph, and even reduce the graph to less than $0.1\%$ of the original size. Besides, we can see that our \textbf{Causal Reduce} achieves best reduction effect in all datasets, that is, our \textbf{Causal Reduce} results in a much smaller kernel size than other algorithms. Moreover, compared with \textbf{Reduce}, \textbf{Re-Confin} can achieve better reduction effect in all datasets, while \textbf{Re-Cover} has basically no performance improvement. This shows that replacing the step~$5$ of \textbf{Reduce} with \textbf{Remove Unconfined \& Contract Confining} plays a key role in improving the performance of the reduction algorithm, and combined with \textbf{Remove Uncovered \& Contract Covering}, the performance of the reduction algorithm will be greatly improved, but only adding \textbf{Remove Uncovered \& Contract Covering} can hardly improve the performance of the reduction algorithm. 

More notably, our \textbf{Causal Reduce} take less time than other algorithms on half of the datasets. On the rest of the datasets, our \textbf{Causal Reduce} only takes a few seconds longer than other algorithms. These phenomena show that integrating our CITs into the reduction algorithm can significantly improve the performance of the algorithm, but the increase in time cost is very small, and they can even reduce the time cost.

\subsection{Performance Gain of CITs on the B\&R Algorithm}
\label{subsec:Performance Analysis of Causal Inference Techniques over Exact Algorithm}
We will examine the performance gain of our CITs for B\&R algorithm. The running time bound is set as $1,000$ seconds for all algorithms, and if the algorithm cannot find the optimal solution within the time bound, the best solution found in all search branches is output.
\begin{table}[ht]
\centering
\resizebox{\textwidth}{27.5mm}{
\begin{tabular}{|c|c|c|c|c|c|c|c|c|c|c|c|c|c|c|c|c|}
\hline 
                   	          &  \multicolumn{2}{c|}{\textbf{Solve}}  &   \multicolumn{2}{c|}{\textbf{Solve-CR}}    & \multicolumn{2}{c|}{\textbf{Solve-Packing}} & \multicolumn{2}{c|}{\textbf{Solve-CR-IC}} & \multicolumn{2}{c|}{\textbf{Causal B\&R Solver}}       \\ \hline 
\textbf{NO.}            & \textbf{Time(S)}& \textbf{Result} &\textbf{Time(S)} &  \textbf{Result}  & \textbf{Time(S)} & \textbf{Result} & \textbf{Time(S)} & \textbf{Result} & \textbf{Time(S)} & \textbf{Result}       \\ \hline
$\mathbf{1}$                  & $1000$          & $1380579565$    &   $1000$        &  $1380810330$     &   $1000$         &  $1380579506$  & $1000$           &   $1380980673$   & $1000$           &   $\bf1380980749$     \\ \hline
$\mathbf{2}$                  & $1000$          & $232813323$     &   $28.5058$     &  $\bf232828253$   &   $1000$         &  $232814520$   & $28.8413$     &   $\bf232828253$  & $\underline{25.4889}$     &   $\bf232828253$      \\ \hline 
$\mathbf{3}$                  & $\underline{3.8108}$    & $\bf10303506$   &   $5.3508$      &  $\bf10303506$    &   $4.7170$      &  $\bf10303506$ & $4.1507$        &   $\bf10303506$       & $3.9148$        &   $\bf10303506$       \\ \hline
$\mathbf{4}$                  & $1000$          & $124020452$     &   $1000$        &  $124020466$      &   $1000$         &  $124020706$    & $1000$           &   $124021474$     & $1000$           &   $\bf124022398$      \\ \hline
$\mathbf{5}$                  & $1000$          & $3904552$       &   $\underline{4.2820}$  &  $\bf3916599$     &   $1000$         &  $3904544$      & $4.3919$        &   $\bf3916599$& $4.33729$        &   $\bf3916599$        \\ \hline
$\mathbf{6}$                  & $1000$          & $100956490$     &   $1000$        &  $101259145$      &   $1000$         &  $100957090$   & $1000$           &   $101259161$ & $1000$           &   $\bf101288073$      \\ \hline
$\mathbf{7}$                  & $1000$          & $798872105$     &   $1000$        &  $799021102$      &   $1000$         &  $798911163$    & $1000$           &   $\bf799021209$& $1000$           &   $\bf799021209$      \\ \hline
$\mathbf{8}$                  & $1000$          & $134613130$     &   $34.0139$     &  $\bf134621271$   &   $1000$         &  $134613130$    & $34.4151$     &   $\bf134621271$ & $\underline{32.1868}$     &   $\bf134621271$      \\ \hline
$\mathbf{9}$                  & $1000$          & $7237240$       &   $\underline{6.8162}$  &  $\bf7273973$     &   $1000$         &  $7237411$     & $6.8372$        &   $\bf7273973$        & $6.9285$        &   $\bf7273973$        \\ \hline
$\mathbf{10}$                 & $1000$          & $71945454$      &   $1000$        &  $71944343$       &   $1000$         &  $\bf71946049$  & $1000$           &   $71944343$  & $1000$           &   $71945241$          \\ \hline
$\mathbf{11}$                 & $1000$          & $2707746$       &   $1000$        &  $\bf2743962$     &   $1000$         &  $2707846$      & $1000$           &   $\bf2743962$ & $1000$           &   $\bf2743962$        \\ \hline
$\mathbf{12}$                 & $1000$          & $61702804$      &   $1000$        &  $61818326$    &   $1000$         &  $61702794$     & $1000$           &   $\bf 61819495$ & $1000$           &   $61818234$       \\ \hline
\end{tabular}
}
\caption{Performance gain of CITs on B\&R algorithm. The \textbf{bold} and \underline{underlined} numbers are the best numerical results of all algorithms and the shortest running time of all algorithms to find the optimal solution, respectively.}
\label{tab:Performance Analysis of Our Causal Inference Techniques. Better results are highlighted in bold.}
\end{table} 

We output the numerical results and running times of all algorithms in Table~\ref{tab:Performance Analysis of Our Causal Inference Techniques. Better results are highlighted in bold.}. It can be seen from Table~\ref{tab:Performance Analysis of Our Causal Inference Techniques. Better results are highlighted in bold.} that \textbf{Solve-CR} and \textbf{Solve-CR-IC}, like our \textbf{Causal B\&R Solver}, can obtain the optimal solution in five data sets, while \textbf{Solve-Packing}, like \textbf{Solve}, can only obtain the optimal solution in one data set. In addition, on those datasets where the optimal solution cannot be solved within $1000$ seconds, our \textbf{Causal B\&R Solver} can basically obtain better numerical solutions than \textbf{Solve-CR-IC}, and \textbf{Solve-CR-IC} can obtain numerical results that are slightly better than \textbf{Solve-CR},  while \textbf{Solve-Packing} can generally get better numerical solutions than \textbf{Solve}. These results demonstrate  that our reduction algorithm, \textbf{Causal Reduce}, is critical for the B\&R algorithm to obtain optimal solutions on more datasets. Moreover, both the inferred covering set of the branching vertex and the weighted packing constraints can help B\&R algorithm find more promising branches and find better solutions.

\subsection{Causal Reduce's Improvement on Heuristic Algorithm}
\label{subsec:Improvement of Causal Reduce for Heuristic Algorithm}
Next, we will verify the superiority of our \textbf{Causal Reduce} as a preprocessing for improving the heuristic algorithm.
Table~\ref{tab:The experimental effect of heuristic algorithm is improved by reduction algorithm. Better results are highlighted in bold.} presents the running time (including preprocessing time) and numerical results.
We find that the preprocessed heuristic algorithm with \textbf{Causal Reduce} usually stop execution after running for a short time, while the rest of the heuristic algorithms are allowed to run for $1000$ seconds. 
Meanwhile, it can be observed from  Table~\ref{tab:The experimental effect of heuristic algorithm is improved by reduction algorithm. Better results are highlighted in bold.} that adding the reduction algorithm as preprocessing is obvious for improving the performance of the heuristic algorithm, and our \textbf{Causal Reduce} helps heuristics find better solutions on all instances in less time (essentially within 100 seconds) than \textbf{Reduce}. Thus, although our \textbf{Causal Reduce} takes no more than $12$ seconds longer than \textbf{Reduce} on half of the datasets~(as can be known from the numerical results in Section~\ref{subsec:Impact of CITs on Reduction Algorithm}), it can further reduce the size of remaining graph by more than $32.6\%$, which is critical for subsequent processing of the problem~(also be mentioned in Section~\ref{subsec:Performance Analysis of Causal Inference Techniques over Exact Algorithm}), so such processing time cost is worth it!

\begin{table}[ht]
\centering
\resizebox{\textwidth}{27.5mm}{
\begin{tabular}{|c|c|c|c|c|c|c||c|c|c|c|c|c|c|c|c|c|c|c|}
\hline 
                & \multicolumn{2}{c|}{\textbf{Fast}}& \multicolumn{2}{c|}{\textbf{Re + Fast}}&\multicolumn{2}{c||}{\textbf{Causal Re + Fast}} & \multicolumn{2}{c|}{\textbf{Dyn}}  & \multicolumn{2}{c|}{\textbf{Re + Dyn}} & \multicolumn{2}{c|}{\textbf{Causal Re + Dyn}}       \\ \hline 
\textbf{NO.}    & \textbf{Time(S)}& \textbf{Result} &\textbf{Time(S)} &  \textbf{Result}     & \textbf{Time(S)} & \textbf{Result}            & \textbf{Time(S)} & \textbf{Result} & \textbf{Time(S)} & \textbf{Result} & \textbf{Time(S)} & \textbf{Result}    \\ \hline
$\mathbf{1}$    &   $1000$        &  $1308864893$   &   $1000$        &  $1381212394$        &  $\underline{250}$       &  $\bf1381215439$            &  $1000$          &  $1310465732$   &   $1000$         &  $1381174854$   & $\underline{100}$         &  $\bf1381178355$   \\ \hline
$\mathbf{2}$    &   $1000$        &  $227121250$    &   $1000$        &  $232826881$         &  $\underline{20}$        &  $\bf232827816$             &  $1000$          &  $229769205$    &   $1000$         &  $232827891$    & $\underline{20}$          &  $\bf232828157$    \\ \hline
$\mathbf{3}$    &   $1000$        &  $10044429$     &   $1000$        &  $10302725$          &  $\underline{5}$         &  $\bf10303465$              &   $1000$         &  $10224463$     &   $1000$         &  $10303168$     & $\underline{5}$           &  $\bf10303476$       \\ \hline
$\mathbf{4}$    &   $1000$        &  $122468973$    &   $1000$        &  $124025600$         &  $\underline{6}$         &  $\bf124026219$             &   $1000$         &  $123179849$    &   $1000$         &  $124026286$    & $\underline{6}$           &  $\bf124026433$      \\ \hline
$\mathbf{5}$    &   $1000$        &  $3823908$      &   $1000$        &  $3916222$           &  $\underline{3}$         &  $\bf3916534$               &   $1000$         &  $3894401$      &   $1000$         &  $3916381$      & $\underline{3}$           &  $\bf3916568$        \\ \hline
$\mathbf{6}$    &   $1000$        &  $97155884$     &   $1000$        &  $101739247$         &  $\underline{275}$       &  $\bf101745579$             &   $1000$         &  $99122831$     &   $1000$         &  $101740261$    & $\underline{250}$         &  $\bf101744463$      \\ \hline
$\mathbf{7}$    &   $1000$        &  $756581612$    &   $1000$        &  $799087671$         & $\underline{150}$        &  $\bf799111311$             &   $1000$         &  $757212666$    &   $1000$         &  $799058579$    & $\underline{70}$          & $\bf799083159$       \\ \hline
$\mathbf{8}$    &   $1000$        &  $131848587$    &   $1000$        &  $134620719$         & $\underline{25}$         &  $\bf134621064$             &   $1000$         &  $132970128$    &   $1000$         &  $134621142$    & $\underline{25}$          & $\bf134621249$       \\ \hline
$\mathbf{9}$    &   $1000$        &  $7131770$      &   $1000$        &  $7273626$           & $\underline{6}$          &  $\bf7273655$               &   $1000$         &  $7252153$      &   $1000$         &  $7273750$      & $\underline{6}$           & $\bf7273931$         \\ \hline
$\mathbf{10}$   &   $1000$        &  $70966490$     &   $1000$        &  $71946839$          & $\underline{6}$          &  $\bf71947488$              &   $1000$         &  $71459210$     &   $1000$         &  $71947516$     & $\underline{6}$           & $\bf71947590$        \\ \hline
$\mathbf{11}$   &   $1000$        &  $2682982$      &   $1000$        &  $2748925$           & $\underline{25}$         &  $\bf2748945$               &   $1000$         &  $2743648$      &   $1000$         &  $2748982$      & $\underline{10}$          & $\bf2749005$         \\ \hline
$\mathbf{12}$   &   $1000$        &  $59619787$     &   $1000$        &  $61850413$          & $\underline{50}$         &  $\bf61852628$              &   $1000$         &  $60802433$     &   $1000$         &  $61855209$     & $\underline{75}$          & $\bf61857313$        \\ \hline
\end{tabular}
}
\caption{A comparative experiment of the effect of \textbf{Causal Reduce} on improving the heuristic algorithm. The numbers in \textbf{bold} and \underline{underlined} are the corresponding experimental results and running time~(including preprocessing time) of each heuristic algorithm using our \textbf{Causal Reduce} as preprocessing, respectively.}
\label{tab:The experimental effect of heuristic algorithm is improved by reduction algorithm. Better results are highlighted in bold.}
\end{table}

\subsection{Comparative Experiment on Causal Search}
\label{subsec:Local search algorithm combined with local causal inference techniques}
On the basis of preprocessing the input graph with \textbf{Causal Reduce}, we will compare our \textbf{Causal Search} with DynWVC2 to verify the effect of adding CITs to the local search process of DynWVC2 algorithm. The running time for both algorithms~(including pre-processing time) is set to $1000$ seconds. To avoid randomness, we run each instance $5$ times and record the mean and maximum values. Furthermore, in order to estimate the gap between the results obtained by these two algorithms and the MWIS, we need to calculate the upper bound of each instance. 
The upper bound for the $2$nd, $3$rd, $5$th, $8$th, $9$th instance is nothing but the weight of the optimal solution obtained by \textbf{Causal B\&R Solver},
and for the rest of the instances, it is obtained by applying the weighted clique cover method mentioned in Section~\ref{subsec:The Causal B&R Solver} to the remaining graph obtained by \textbf{Causal Reduce}.  
Table~\ref{tab:Local search algorithm combined with causal inference technique. Better results are highlighted in bold.} outputs the numerical results and the estimated gap. The small gaps there demonstrate that after preprocessing with our \textbf{Causal Reduce}, both algorithms can obtain numerical results very close to the optimal solution. In particular, for those instances where the optimal solution is obtained, their gap can basically reach $10^{-6} \sim 10^{-7}$, and in the remaining instances, the estimated gap can basically reach $10^{-4} \sim 10^{-2}$. Besides, from the mean and maximum values,  our \textbf{Causal Search} can basically achieve better performance than DynWVC2, thereby implying that our CITs 
can help local search find better local optima.


\begin{table}[h]
\centering
\small
\resizebox{\textwidth}{25mm}{
\begin{tabular}{|c|c|c|c|c|c|c|c|c|c|c|c|c|c|c|c|c|c|}
\hline 
                  &                      &    \multicolumn{4}{c|}{\textbf{Dyn}}                                                      &   \multicolumn{4}{c|}{\textbf{Causal Search}}                                        \\ \hline 
\textbf{NO.}      & \textbf{Upper Bound} &     \textbf{Mean}     &  \textbf{Gap}   &    \textbf{Max} &  \textbf{Gap}   &   \textbf{Mean}                &  \textbf{Gap}   & \textbf{Max}      &  \textbf{Gap}         \\ \hline
$\mathbf{1}$      &   $1384376268$       &    $1381464698.6$     & $2.103 \times 10^{-3}$ &  $1381467306$   & $2.101 \times 10^{-3}$ &   $\underline{1381466033.6}$   & $2.102 \times 10^{-3}$ & $\bf1381470750$   & $2.099 \times 10^{-3}$ \\ \hline
$\mathbf{2}$      &   $232828253$        &    $232828153.4$      & $4.278 \times 10^{-7}$ &  $232828171$    & $3.522 \times 10^{-7}$ &   $\underline{232828159.2}$    & $4.029 \times 10^{-7}$ & $\bf232828188$    & $2.792 \times 10^{-7}$ \\ \hline
$\mathbf{3}$      &    $10303506$        &    $10303485.2$       & $2.019 \times 10^{-6}$ &  $10303491$     & $1.456 \times 10^{-6}$ &   $\underline{10303488.6}$     & $1.689 \times 10^{-6}$ & $\bf10303494$     & $1.165 \times 10^{-6}$ \\ \hline
$\mathbf{4}$      &    $124076790$       &    $124026438.6$      & $4.058 \times 10^{-4}$ &  $\bf124026451$ & $4.057 \times 10^{-4}$ &   $\underline{124026444.8}$    & $4.058 \times 10^{-4}$ & $124026449$       & $4.057 \times 10^{-4}$ \\ \hline
$\mathbf{5}$      &    $3916599$         &    $3916582.2$        & $4.289 \times 10^{-6}$ &  $3916583$      & $4.085 \times 10^{-6}$ &   $\underline{3916582.6}$      & $4.187 \times 10^{-6}$ & $\bf3916584$      & $3.830 \times 10^{-6}$ \\ \hline
$\mathbf{6}$      &    $103562461$       &    $101846521.0$      & $1.657 \times 10^{-2}$  &  $101848242$    & $1.655 \times 10^{-2}$  &   $\underline{101847524.8}$    & $1.656 \times 10^{-2}$ & $\bf101849650$    & $1.654 \times 10^{-2}$ \\ \hline
$\mathbf{7}$      &   $800748442$        &    $799264479.2$      & $1.853 \times 10^{-3}$ &  $799265827$    & $1.852 \times 10^{-3}$ &   $\underline{799266278.8}$    & $1.851 \times 10^{-3}$ & $\bf799267573$    & $1.849 \times 10^{-3}$ \\ \hline
$\mathbf{8}$      &   $134621271$        &    $134621255.0$      & $1.189 \times 10^{-7}$ &  $134621257$    & $1.040 \times 10^{-7}$ &   $\underline{134621256.8}$    & $1.055 \times 10^{-7}$ & $\bf134621265$    & $4.457\times 10^{-6}$ \\ \hline
$\mathbf{9}$      &    $7273973$         &$\underline{7273939.8}$& $4.564 \times 10^{-6}$ &  $7273945$      & $3.849 \times 10^{-6}$ &   $7273936.8$                  & $4.977 \times 10^{-6}$ & $\bf7273947$      & $3.574 \times 10^{-6}$ \\ \hline
$\mathbf{10}$     &    $71978922$        &    $71947636.0$       & $4.347 \times 10^{-4}$ &  $71947639$     & $4.346 \times 10^{-4}$ &   $\underline{71947636.8}$     & $4.346 \times 10^{-4}$ & $\bf71947642$     & $4.346 \times 10^{-4}$  \\ \hline
$\mathbf{11}$     &    $2819343$         &    $2749008.4$        & $2.495 \times 10^{-2}$  &  $2749009$      & $2.495 \times 10^{-2}$  &   $\underline{2749010.8}$      & $2.495 \times 10^{-2}$ & $\bf2749019$      & $2.494 \times 10^{-2}$ \\ \hline
$\mathbf{12}$     &    $62276875$        &    $61860510.4$       & $6.686 \times 10^{-3}$ & $61860616$      & $6.684 \times 10^{-3}$ &   $\underline{61860549.0}$     & $6.685 \times 10^{-3}$ & $\bf61860717$     & $6.682 \times 10^{-3}$ \\ \hline
\end{tabular}
}
\caption{Compare our \textbf{Causal Search} with the DynWVC2 algorithm. The \textbf{bold} and \underline{underlined} numbers are better maximum and average values, respectively. In the column headed by ``\textbf{Upper Bound}", each number is the upper bound of the MWIS of the corresponding instance.
}
\label{tab:Local search algorithm combined with causal inference technique. Better results are highlighted in bold.}
\end{table}

\section{Conclusion and Outlook}
\label{sec:Conclusion}

In this paper, we propose a series of causal inference techniques~(CITs) for the maximum weight independent set~(MWIS) problem by fully exploiting the upper bound property of MWIS. After integrating our CITs, the performance of various existing algorithms, including the Branch-and-Reduce~(B\&R) algorithm and some heuristic algorithms, is significantly improved. We are now conducting theoretical analysis to find some guarantees on solution quality, developing strategies to help the B\&R algorithm analyze the causes of conflicts and perform more efficient backtracking searches, and generalizing the proposed CITs to other combinatorial optimization problems.

\section*{Acknowledgements}
This research was supported by the National Key R\&D Program of China (Nos. 2020AAA0105200, 2022YFA1005102) and the National Natural Science Foundation of China (Nos.~12288101, 11822102).
SS is partially supported by Beijing Academy of Artificial Intelligence (BAAI). The authors would like to thank Professor Hao Wu for his useful discussions and valuable suggestions.

\normalem

\begin{thebibliography}{10}

\bibitem{akiba2016branch}
{\sc T.~Akiba and Y.~Iwata}, {\em Branch-and-reduce exponential/fpt algorithms
  in practice: A case study of vertex cover}, Theoretical Computer Science, 609
  (2016), pp.~211--225.

\bibitem{avenali2007resolution}
{\sc A.~Avenali}, {\em Resolution branch and bound and an application: the
  maximum weighted stable set problem}, Operations research, 55 (2007),
  pp.~932--948.

\bibitem{babel1994fast}
{\sc L.~Babel}, {\em A fast algorithm for the maximum weight clique problem},
  Computing, 52 (1994), pp.~31--38.

\bibitem{balas1986finding}
{\sc E.~Balas and C.~S. Yu}, {\em Finding a maximum clique in an arbitrary
  graph}, SIAM Journal on Computing, 15 (1986), pp.~1054--1068.

\bibitem{cai2015balance}
{\sc S.~Cai}, {\em Balance between complexity and quality: Local search for
  minimum vertex cover in massive graphs}, in Twenty-Fourth International Joint
  Conference on Artificial Intelligence, 2015, pp.~747--753.

\bibitem{cai2018improving}
{\sc S.~Cai, W.~Hou, J.~Lin, and Y.~Li}, {\em Improving local search for
  minimum weight vertex cover by dynamic strategies.}, in Twenty-Seventh
  International Joint Conference on Artificial Intelligence, 2018,
  pp.~1412--1418.

\bibitem{cai2019towards}
{\sc S.~Cai, Y.~Li, W.~Hou, and H.~Wang}, {\em Towards faster local search for
  minimum weight vertex cover on massive graphs}, Information Sciences, 471
  (2019), pp.~64--79.

\bibitem{chang2017computing}
{\sc L.~Chang, W.~Li, and W.~Zhang}, {\em Computing a near-maximum independent
  set in linear time by reducing-peeling}, in Proceedings of the 2017 ACM
  International Conference on Management of Data, 2017, pp.~1181--1196.

\bibitem{coniglio2022optimizing}
{\sc S.~Coniglio and S.~Gualandi}, {\em Optimizing over the closure of rank
  inequalities with a small right-hand side for the maximum stable set problem
  via bilevel programming}, INFORMS Journal on Computing, 34 (2022),
  pp.~1006--1023.

\bibitem{dahlum2016accelerating}
{\sc J.~Dahlum, S.~Lamm, P.~Sanders, C.~Schulz, D.~Strash, and R.~F. Werneck},
  {\em Accelerating local search for the maximum independent set problem}, in
  International symposium on experimental algorithms, 2016, pp.~118--133.

\bibitem{dzulfikar2019pace}
{\sc M.~A. Dzulfikar, J.~K. Fichte, and M.~Hecher}, {\em The pace 2019
  parameterized algorithms and computational experiments challenge: The fourth
  iteration}, in 14th International Symposium on Parameterized and Exact
  Computation, vol.~148, 2019, pp.~25:1--25:23.

\bibitem{feo1994greedy}
{\sc T.~A. Feo, M.~G. Resende, and S.~H. Smith}, {\em A greedy randomized
  adaptive search procedure for maximum independent set}, Operations Research,
  42 (1994), pp.~860--878.

\bibitem{fomin2009measure}
{\sc F.~V. Fomin, F.~Grandoni, and D.~Kratsch}, {\em A measure \& conquer
  approach for the analysis of exact algorithms}, Journal of the ACM, 56
  (2009), pp.~1--32.

\bibitem{kneis2009fine}
{\sc J.~Kneis, A.~Langer, and P.~Rossmanith}, {\em A fine-grained analysis of a
  simple independent set algorithm}, in IARCS Annual Conference on Foundations
  of Software Technology and Theoretical Computer Science, vol.~4, 2009,
  pp.~287--298.

\bibitem{lamm2019exactly}
{\sc S.~Lamm, C.~Schulz, D.~Strash, R.~Williger, and H.~Zhang}, {\em Exactly
  solving the maximum weight independent set problem on large real-world
  graphs}, in Proceedings of the Twenty-First Workshop on Algorithm Engineering
  and Experiments, 2019, pp.~144--158.

\bibitem{li2018incremental}
{\sc C.-M. Li, Z.~Fang, H.~Jiang, and K.~Xu}, {\em Incremental upper bound for
  the maximum clique problem}, INFORMS Journal on Computing, 30 (2018),
  pp.~137--153.

\bibitem{liao2016approximation}
{\sc C.-S. Liao, C.-W. Liang, and S.~H. Poon}, {\em Approximation algorithms on
  consistent dynamic map labeling}, Theoretical Computer Science, 640 (2016),
  pp.~84--93.

\bibitem{nogueira2018hybrid}
{\sc B.~Nogueira, R.~G. Pinheiro, and A.~Subramanian}, {\em A hybrid iterated
  local search heuristic for the maximum weight independent set problem},
  Optimization Letters, 12 (2018), pp.~567--583.

\bibitem{plachetta2021sat}
{\sc R.~Plachetta and A.~van~der Grinten}, {\em Sat-and-reduce for vertex
  cover: Accelerating branch-and-reduce by sat solving}, in Proceedings of the
  Workshop on Algorithm Engineering and Experiments, 2021, pp.~169--180.

\bibitem{pullan2009optimisation}
{\sc W.~Pullan}, {\em Optimisation of unweighted/weighted maximum independent
  sets and minimum vertex covers}, Discrete Optimization, 6 (2009),
  pp.~214--219.

\bibitem{rossi2015network}
{\sc R.~Rossi and N.~Ahmed}, {\em The network data repository with interactive
  graph analytics and visualization}, in Twenty-ninth AAAI conference on
  artificial intelligence, 2015, pp.~4292--4293.

\bibitem{san2019new}
{\sc P.~San~Segundo, F.~Furini, and J.~Artieda}, {\em A new branch-and-bound
  algorithm for the maximum weighted clique problem}, Computers $\&$ Operations
  Research, 110 (2019), pp.~18--33.

\bibitem{sewell1998branch}
{\sc E.~C. Sewell}, {\em A branch and bound algorithm for the stability number
  of a sparse graph}, INFORMS Journal on Computing, 10 (1998), pp.~438--447.

\bibitem{shyu2004ant}
{\sc S.~J. Shyu, P.-Y. Yin, and B.~M. Lin}, {\em An ant colony optimization
  algorithm for the minimum weight vertex cover problem}, Annals of Operations
  Research, 131 (2004), pp.~283--304.

\bibitem{strash2016power}
{\sc D.~Strash}, {\em On the power of simple reductions for the maximum
  independent set problem}, in International Computing and Combinatorics
  Conference, vol.~9797, Springer, 2016, pp.~345--356.

\bibitem{warren2006combinatorial}
{\sc J.~S. Warren and I.~V. Hicks}, {\em Combinatorial branch-and-bound for the
  maximum weight independent set problem}, Relat{\'o}rio T{\'e}cnico, Texas
  A$\&$M University, Citeseer, 9 (2006), p.~17.

\bibitem{xiao2021efficient}
{\sc M.~Xiao, S.~Huang, Y.~Zhou, and B.~Ding}, {\em Efficient reductions and a
  fast algorithm of maximum weighted independent set}, in Proceedings of the
  Web Conference 2021, 2021, pp.~3930--3940.

\bibitem{xiao2013confining}
{\sc M.~Xiao and H.~Nagamochi}, {\em Confining sets and avoiding bottleneck
  cases: A simple maximum independent set algorithm in degree-3 graphs},
  Theoretical Computer Science, 469 (2013), pp.~92--104.

\bibitem{xu2016new}
{\sc H.~Xu, T.~Kumar, and S.~Koenig}, {\em A new solver for the minimum
  weighted vertex cover problem}, in International Conference on AI and OR
  Techniques in Constraint Programming for Combinatorial Optimization Problems,
  vol.~9676, 2016, pp.~392--405.

\end{thebibliography}

\end{document}